\documentclass[12pt]{amsart}
\usepackage{pb-diagram}
\usepackage{lamsarrow}
\usepackage{pb-lams}
\usepackage{euscript}
\usepackage{cmex12-2e}
\usepackage[6in]{amspapermacs}

\def\modulefont#1{\mathsf{#1}}
\def\modulebundfont#1{\bundlefont{#1}}
\def\X{{\modulefont X}}
\def\bX{{\modulebundfont X}}
\def\Y{{\modulefont Y}}
\def\bY{{\modulebundfont Y}}
\def\ZZ{{\modulefont Z}}
\def\bZZ{{\modulebundfont Z}}
\def\NN{{\modulefont N}}
\def\bNN{{\modulebundfont N}}
\def\smeover #1{\,\mathord{\mathop{\text{--}}\nolimits_{#1}}\,}
\def\sme{\,\mathord{\mathop{\text{--}}\nolimits_{\relax}}\,}
\def\ib{im\-prim\-i\-tiv\-ity bi\-mod\-u\-le}
\def\ibb{\ib{} bundle}
\def\vg{V_\gamma}

\def\bundlefont#1{\EuScript{#1}}
\def\A{{\bundlefont A}}
\def\B{{\bundlefont B}}
\def\CC{{\bundlefont C}}
\def\D{{\bundlefont D}}
\def\TT{{\bundlefont T}}
\def\pa{p_{\A}}
\def\pb{p_{\B}}
\def\pz{p_{\bZZ}}
\def\ax{\A^X}
\def\px{p^X}
\def\bA{\overline{\A}}
\def\bp{\bar p}
\DeclareMathOperator{\br}{Br}
\DeclareMathOperator{\sbr}{\mathfrak{Br}}
\def\sbrg{\sbr(G)}
\def\brg{\br(G)}
\def\sbrz{\sbr_0}
\def\brz{\br_0}
\def\brzg{\brz(G)}
\def\sbrzg{\sbrz(G)}
\def\aa{(\A,\a)}
\def\bb{(\B,\beta)}
\def\xv{{(\bX,V)}}
\def\oX{\overline{\X}}
\def\K{K}
\def\KK{{\bundlefont K}}
\def\HH{{\bundlefont H}}
\def\H{\modulefont{H}}
\def\hs{H}
\def\kh{\K(\H)}
\def\kkhh{\KK(\HH)}

\def\unitspace#1{{{#1}^{(0)}}}
\def\go{\unitspace G}
\def\ho{\unitspace H}
\def\eo{\unitspace E}
\def\hmgrp{H\backslash\grp}
\def\grpmg{\grp/G}
\def\g{\gamma}
\def\gi{{\g_i}}
\def\rx{r_X}
\def\sx{s_X}
\def\equi(#1,#2){$(#1,#2)$-equi\-va\-lence}
\def\opposite#1{{#1}^{\text{op}}}
\def\xop{\opposite{X}}
\def\imp(#1){#1\understar{G}\opposite{#1}}
\def\xsgy{{X\understar{G}Y}}
\def\xopshx{\xop\understar HX}
\def\xsgxop{X\understar G\xop}
\DeclareMathOperator{\Ext}{Ext}
\DeclareMathOperator{\E}{\mathcal{E}}
\DeclareMathOperator{\PP}{\mathcal{P}}
\def\pg{\PP(G)}
\def\lbb{\lbrack\!\lbrack}
\def\rbb{\rbrack\!\rbrack}
\def\eclass#1{\lbb #1 \rbb}
\def\ext(#1,#2){\Ext(#1,#2)}
\def\extgt{\ext(G,\T)}
\DeclareMathOperator{\Gr}{Gr}
\def\grp{{\Gr(\phi)}}
\DeclareMathOperator{\Tw}{Tw}
\def\myatop#1#2#3{\mathrel{\genfrac{}{}{0pt}{#3}{#2}{#1}}}
\def\understar#1{\mathrel{\mathchoice{%
\myatop{#1}{\textstyle{*}}1}%
{\myatop{#1}{\textstyle{*}}1}%
{\myatop{#1}{\scriptstyle{*}}2}%
{\myatop{#1}{\scriptscriptstyle{*}}2}%
}}
\def\stt{\understar{\T}}
\def\ste{\understar{E}}
\def\stg{\understar{G}}
\def\teq{$\T$-equi\-va\-lence}
\def\simh{\sim_H}
\def\efteq{$(E,F)$-\teq}
\def\tphi{\tilde\phi}
\def\ttphi{\Tilde{\Tilde{\phi}}}
\def\ttpsi{\Tilde{\Tilde{\psi}}}
\def\eq[#1,#2]{\bigl[#1,\iota(#2)\bigr]}
\def\smeq[#1,#2]{[#1,\iota(#2)]}

\def\ag{\alpha_\gamma}
\def\a{\alpha}
\def\ba{\bar\a}

\def\batgob{\A\tensor_{\go} \B}
\def\atgob{A\tensor_{\go} B}

\def\U{{\boldsymbol U}}
\def\Fi{{F_i}}
\def\fij{{F_{ij}}}
\def\gij{{g_{ij}}}
\def\L{\mathcal{L}}
\def\sT{\mathcal{S}}
\def\tb{\Lambda}
\def\btb{\bar\tb}
\def\tbstbtb{\tb\stt\btb}

\makeatletter
\@addtoreset{footnote}{section}
\makeatother

\makeatletter
\def\labelenumi{(\@alph\c@enumi)}
\def\theenumi{\@alph \c@enumi}
\newcount\charno
\def\alphapart#1{\charno=96
\advance\charno by#1\char\charno}
\def\partref#1{\textup{(}\alphapart{#1}\textup{)}}
\makeatother
\def\<{\langle}
\def\>{\rangle}
\let\ipscriptstyle=\scriptscriptstyle
\def\lipsqueeze{{\mskip -3.0mu}}
\def\ripsqueeze{{\mskip -3.0mu}}
\def\ipcomma{\mathrel{,}}
\newbox\ipstrutbox
\setbox\ipstrutbox=\hbox{\vrule height8.5pt
width 0pt}
\def\ipstrut{\copy\ipstrutbox}
\def\lip#1<#2,#3>{\mathopen{\relax_{\ipstrut\ipscriptstyle{
#1}}\lipsqueeze
\langle} #2\ipcomma #3 \rangle}
\def\blip#1<#2,#3>{\mathopen{\relax_{\ipstrut
\ipscriptstyle{ #1}}\lipsqueeze\bigl\langle} #2\ipcomma #3 \bigr\rangle}
\def\rip#1<#2,#3>{\langle #2\ipcomma #3
\rangle_{\ripsqueeze\ipstrut\ipscriptstyle{#1}}}
\def\brip#1<#2,#3>{\bigl\langle #2\ipcomma #3
\bigr\rangle_{\ripsqueeze\ipstrut\ipscriptstyle{#1}}}
\def\angsqueeze{\mskip -6mu}
\def\smangsqueeze{\mskip -3.7mu}
\def\trip#1<#2,#3>{\langle\smangsqueeze\langle #2\ipcomma #3
\rangle\smangsqueeze\rangle_{\ripsqueeze\ipstrut\ipscriptstyle{#1}}}
\def\btrip#1<#2,#3>{\bigl\langle\angsqueeze\bigl\langle #2\ipcomma
#3
\bigr\rangle
\angsqueeze\bigr\rangle_{\ripsqueeze\ipstrut\ipscriptstyle{#1}}}
\def\tlip#1<#2,#3>{\mathopen{\relax_{\ipstrut\ipscriptstyle{
#1}}\lipsqueeze \langle\smangsqueeze\langle} #2\ipcomma #3
\rangle\smangsqueeze\rangle}
\def\btlip#1<#2,#3>{\mathopen{\relax_{\ipstrut\ipscriptstyle{
#1}}\lipsqueeze
\bigl\langle\angsqueeze\bigl\langle} #2\ipcomma #3
\bigr\rangle\angsqueeze\bigr\rangle}

\def\ip(#1|#2){(#1\mid #2)}
\def\bip(#1|#2){\(#1 \mid #2\)}

\def\psnrl{pass to a subnet and relabel}
\def\hm{homomorphism}
\def\sclcg{second countable locally compact groupoid}
\def\eltwo{\ell^2}
\def\rg{\mathcal{R}(G)}
\def\simext{\sim_{\text{ext}}}
\def\simcoh{\sim_{\text{coh}}}
\def\Th#1{\Theta_{#1}}
\def\thetag{\Th{G}}
\def\l{\lambda}
\def\pil{\pi_\ell}

\def\gotoco(#1,#2){C_0(#1;#2)} 
\def\gsuboT(#1){\gotoco(T,#1)}
\def\gsubogo(#1){\gotoco(\go,#1)}
\def\spull #1^#2{{#1}^{#2}}

\begin{document}

\title[The Brauer Group]{The Brauer Group of a Locally Compact Groupoid}

\author[Kumjian]{Alexander Kumjian}
\address{Department of Mathematics \\
University of Nevada \\
Reno, NV 89557 \\
USA}
\email{alex@math.unr.edu}

\author[Muhly]{Paul S. Muhly}
\address{%
Department of Mathematics \\
University of Iowa \\
Iowa City, Iowa 52242}

\email{muhly@math.uiowa.edu}

\author[Renault]{Jean N. Renault}

\address{D\'epartement de Math\'ematiques \\
Universit\'e d'Orl\'eans \\ 45067 Orl\'eans Cedex 2 \\ France}

\email{renault@univ-orleans.fr}

\author[Williams]{\\Dana P. Williams}
\address{Department of Mathematics \\
Dartmouth College \\
Hanover, New Hampshire 03755-3551 \\
USA}
\email{dana.williams@dartmouth.edu}

\thanks{The second author was partially supported by
the National Science Foundation.}

\date{20 January 1997}



\begin{abstract}
%
%

We define the Brauer group $\brg$ of a locally compact groupoid $G$ to
be the set of Morita equivalence classes of pairs $\aa$ consisting of
an elementary \cs-bundle $\A$ over $\go$ satisfying Fell's condition
and an action $\a$ of $G$ on
$\A$ by $*$-isomorphisms.  When $G$ is the transformation groupoid
$X\times H$, then $\brg$ is the equivariant Brauer group $\br_H(X)$.  

In addition to proving that $\brg$ is a group, we prove three
isomorphism results.  First we show that if $G$ and $H$ are
equivalent groupoids, then $\brg$ and $\br(H)$ are isomorphic.  This
generalizes the result that if $G$ and $H$ are groups acting freely
and properly on a space $X$, say $G$ on the left and $H$ on the right
then $\br_G(X/H)$ and $\br_H(G\backslash X)$ are isomorphic.  Secondly
we show that the subgroup $\brzg$ of $\brg$ consisting of classes
$[\A,\a]$ with $\A$ having trivial Dixmier-Douady invariant is
isomorphic to a quotient $\E(G)$ of the collection $\Tw(G)$ of twists
over $G$.  Finally we prove that $\brg$ is isomorphic to the inductive
limit $\extgt$ of the groups $\E(G^X)$ where $X$ varies over all
principal $G$ spaces $X$ and $G^X$ is the imprimitivity groupoid
associated to $X$.

\end{abstract}

\maketitle

\ifdraft
%
%

\section{Introduction}

This paper is about two groups that are naturally associated to a
locally compact groupoid (with Haar system) and an isomorphism between
them. The first group, denoted $\brg$ and called the \emph{Brauer
group of} $G$, where $G$ is the groupoid in question, is a collection
of equivalence classes of actions of $G$ on certain bundles of
$\cs$-algebras. The second group, $\extgt$, is a group that organizes all
the extensions by $\T$ of all the groupoids that are equivalent to $G$
in the sense of
\cite{mrw}. The isomorphism between these two groups or, more accurately, our
analysis of it, may be viewed as a simultaneous generalization of Mackey's
analysis of projective representations of locally compact groups \cite{m4}%
, on the one hand, and the Dixmier-Douady analysis of continuous trace $%
\cs $-algebras (\cite{dd} and \cite{dixfields}), on the other. Moreover, our
work provides a global perspective from which to view a number of recent
studies in the theory of $\cs$-dynamical systems. See, in particular,
\cite{ckrw}, \cite{krw}, \cite{masuda}, \cite{prw}, \cite{rr}, 
\cite{rt}, \cite{rw}, and \cite{90c}.
Further, it suggests a way to generalize to groupoids Moore's Borel
cohomology theory for locally compact groups without becoming
entangled in measure-theoretic difficulties.

To expand upon these remarks, and to help motivate further the need
for our theory, we shall begin with Mackey's analysis \cite{m4}, but
first, we wish to emphasize the usual separability assumptions that
will be made throughout this paper. All groups, groupoids, and spaces
will be locally compact, Hausdorff, and second countable. Thus, in
particular, they are paracompact.  All Hilbert spaces will be
separable. They may be finite dimensional as well as infinite
dimensional. Likewise, all $\cs$-algebras under discussion will be
assumed to be separable.
Since we are interested mainly in groupoids with Haar systems, we
always assume our groupoids have open range and source maps.

Mackey was inspired, in part, 
by earlier investigations of Wigner \cite{ewigner}
and Bargmann \cite{v-barg} concerned with the mathematical foundations of
quantum mechanics and, in particular, with the problem of classifying the
so-called \emph{automorphic representations} of locally compact groups on
(infinite dimensional) Hilbert space. The problem (formulated a bit
differently than in \cite{ewigner}, \cite{v-barg}, and \cite{m4}) 
is to classify the continuous
representations $\alpha $ of a locally compact group
$G$ as automorphisms of the $%
\cs$-algebra $\K=\K(\hs)$ of all compact operators on a separable infinite
dimensional Hilbert space $\hs$. The meaning of `continuity' is that for each
compact operator $k\in \K$, the function $g\rightarrow \alpha _g(k)$ is
continuous from $G$ to $\K$ with the norm topology. Two such representations $%
\alpha _1$ and $\alpha _2$ are considered equivalent if there is a $*$%
-automorphism $\varphi :\K\rightarrow \K$ such that $\varphi \circ \alpha
_1=\alpha _2\circ \varphi $. The collection of equivalence classes of such
representations will be denoted $\mathcal{R}(G)$. With respect to the
process of forming tensor products, it is not difficult to see that $%
\mathcal{R}(G)$ becomes a semigroup with identity
in a natural way. The problem is to describe
this semigroup.

Owing to Wigner's theorem \cite{ewigner} that all $*$-automorphisms of $\K(\hs)$
are implemented by unitary operators on the underlying Hilbert space $\hs$,
and the fact that two unitaries on $\hs$ implement the same automorphism if
and only if one is a unimodular scalar multiple of the other, one may
summarize part of the analysis as follows: Every $\alpha $ is determined by
a unimodular \emph{Borel} function $\sigma $ on $G\times G$ satisfying
identities making it a $2$-cocycle in a certain cohomology theory and a
unitary $\sigma $-representation $U$ of $G$ with values in the unitary group
of $\hs$. (These are also called \emph{projective representations} of $G$.)
The relation between $U$ and $\sigma $ is found in the equation $%
U_{gh}=\sigma (g,h)U_gU_h$, for $g,h\in G$. The formula for 
$\alpha $, then, is $%
\alpha _g(k)=U_gkU_g^{-1}$, $k\in K$, for $g\in G$. Two $\alpha $'s are
equivalent if and only if they are spatially equivalent, and in terms of the 
$U$'s and $\sigma $'s, this happens if and only if $(U_1,\sigma _1)$ is
unitarily equivalent (in the obvious sense) to $(U_2,\sigma _2^{\prime })$
where $\sigma _2^{\prime }$ is \emph{cohomologous} to $\sigma _2$, meaning
that there is a unimodular \emph{Borel} function $b$ on $G$ with the
property that $\sigma _2^{\prime }(g,h)=b(gh)\overline{b(g)b(h)}\cdot
\sigma _2(g,h)$. 
The upshot is that to understand $\mathcal{R}(G)$, one must identify $H^2(G,%
\T )$, and for each $[\sigma ]\in H^2(G,\T )$, one must identify the
unitary equivalence classes of $\sigma $-unitary representations of $G$.

As we are emphasizing, the group $H^2(G,\T )$ is a measure-theoretic
object that cannot be avoided, if one wants to think in terms of projective
representations of $G$. As Mackey noted in \cite{m4}, however, one may
avoid measure theoretic difficulties (or at least push them into the
background), if one thinks in terms of \emph{extensions} and an $\Ext$-group,
instead of $H^2(G,\T )$. The point is, thanks to Mackey 
\cite[Theorem 7.1]{m5} and Weil \cite{weil}, the correspondence between
extensions and $2$-cohomology that one learns about in finite group theory
extends to the theory of locally compact groups as follows. With each
two-cocycle $\sigma $ one can form a certain \emph{locally compact
topological} group $G^\sigma $. Set-theoretically $G^\sigma $ is the product
of the circle $\T $ with $G$, $\T \times G$, but with the
multiplication $(t_1,g_1)(t_2,g_2)=(t_1t_2\sigma (g_1,g_2),g_1g_2)$. In
general, the topology on $G^\sigma $ is \emph{not} the product topology. One
focuses on the Borel structure of $G^\sigma $, i.e., the product Borel
structure, and notes that Haar measure on $\T $ times Haar measure on $G$
endows $G^\sigma $ with a left invariant Borel measure. The theorems of
Mackey and Weil, then, force the existence of the topology. One can prove
without difficulty that two cocycles give rise to topologically isomorphic
extensions if and only if they are cohomologous. Thus, one obtains, as
anticipated, a natural isomorphism between $\extgt$ and $H^2(G,\T)$.

The measure theoretic technology is pushed further into the background when
one recognizes how to realize directly an extension, given an automorphic
representation $\alpha $ of $G$. Indeed, if $\hs$ is the underlying Hilbert
space and if $\mathcal{U}(\hs)$ denotes the set of unitary operators on $\hs$
endowed with the $*$-strong operator topology, simply form $E(\alpha
):=\{(g,U)\in G\times \mathcal{U}(\hs)|\Ad(U)=\alpha _g\}$. It is not hard to
see that with the relative topology, $E(\alpha )$ is a locally compact group
and an extension of $G$ by $\T $. Further, the relation between $%
E(\alpha )$ and $G^\sigma $ is given, via the $\sigma $-unitary
representation, $U^\sigma $, of $G$ determined by $\alpha $, by the formula 
\begin{equation*}
(t,g)\rightarrow (g,tU_g^\sigma ),\quad\text{$(t,g)\in G^\sigma $.}
\end{equation*}

It is not difficult to see that two elements $\alpha ,\beta \in \mathcal{R}%
(G)$ give isomorphic extensions $E(\alpha )$ and $E(\beta )$ if and only if $%
\alpha $ and $\beta $ are \emph{exterior equivalent} in the sense that there
is a unitary valued function $u$ on $G$ such that $\beta _g=Ad(u_g)\circ
\alpha _g$ for all $g\in G$ and such that $u_{gh}=u_g\alpha _g(u_h)$ for all 
$g,h\in G$. Further, it is not difficult to see that $\alpha $ is exterior
equivalent to the identity (i.e., the trivial action) if and only if $\alpha 
$ is implemented by a \emph{unitary representation} of $G$. Thus, the map $%
\alpha \rightarrow E(\alpha )$ sets up an isomorphism between $\mathcal{R}%
(G)/\mathcal{S}$, where $\mathcal{S}$ is the subsemigroup of spatially
implemented automorphic representations, and $\extgt$. This shows, in
particular, that $ {\mathcal R}(G)/{\mathcal S}$ has the structure of a group.
Further, the class of $E(\alpha )$ in $\extgt$ (and also the class of
the cocycle $\sigma $ corresponding to $\alpha $) has come to be known
as the \emph{Mackey obstruction} to implementing $\alpha $ by a
unitary representation of $G$. For reasons that will be further
clarified in a moment, we shall write $%
\brg$ for $\mathcal{R}(G)/\mathcal{S}$ and call this group 
the \emph{Brauer group
of $G$}. We thus have the isomorphisms 
\begin{equation*}
\brg\cong H^2(G,\T )\cong \extgt.
\end{equation*}

We turn now to the second source of motivation for this paper, the
Dixmier-Douady theory of continuous trace $\cs$-algebras (\cite{dd}
and \cite {dixfields}, these are summarized in Chapter~10 of
\cite{dix}). However, we shall 
adopt the somewhat more contemporary perspective first articulated by Green
in \cite{green1} and by Taylor \cite{joe}. (What follows, is amplified in
\cite{ckrw}, which is based in turn on \cite{90c}.) The quickest way to say that a $%
\cs$-algebra $A$, say, is continuous trace is to say that it has Hausdorff
spectrum $T$ and that there is a bundle $\mathcal{A}$ of $\cs$-algebras
over $T$, each fibre $A(t)$ of which is isomorphic to the algebra of compact
operators on some Hilbert space\footnote{%
Such bundles are called \emph{elementary} $\cs$-bundles. It should be kept
in mind that fibres can be finite dimensional as well as infinite
dimensional.}, and satisfying \emph{Fell's condition} (see
Definition~\ref{def-2.11}, 
below), such that $A=\gotoco(T,\A)$, the $\cs$-algebra of continuous
sections of $\A$ that vanish at infinity on $T$. 
If $A$ is unital, then $T$ is compact
and all the fibres of $\A$ are unital (and therefore finite
dimensional). Evidently, in this
case, the center of $A$ is $C(T)$. 
Further in this case, as Grothendieck observed \cite{groth2}, $A$
may be viewed as a separable algebra over $C(T)$ in the sense of pure
algebra. One is then led to think about classifying algebras like $A$ up to
Morita equivalence over $C(T)$. When one does this, the result has a group
structure, denoted $\br(T)$, that generalizes the \emph{Brauer group} of
central simple algebras over a field. Green recognized in \cite{green4} that
this idea works in the non-unital setting as well, provided one uses
Rieffel's notion of strong Morita equivalence \cite{rieff2}. 

In a bit more detail, let us consider continuous trace $\cs$-algebras $A$
with a prescribed spectrum $T$. Then we shall say that two such algebras
$A_1$ and $A_2$ are strongly Morita equivalent (in the sense of Rieffel
\cite{rieff2}) via an equivalence $\X$ \emph{that respects} $T$, or more
simply that $\X$ is a \emph{strong Morita equivalence over} $T$, if $%
\X$ implements a strong Morita equivalence between $A_1$ and $A_2$
in such a way that the actions of $C_0(T)$ on $\X$ are the
same, where $C_0(T)$ is viewed as the center of each of $A_1$ and $A_2$. If $%
A_i$ is realized as $\gotoco(T,\A_i)$, $i=1,2$, then such an $\X$
exists if and only if there is a bundle $\bX$ over $T$ such that for
each $t\in T$, $\X(t)$ implements a strong Morita equivalence
between $A_1(t)$ and $A_2(t)$ and such that $\X=\gotoco(T,\bX)$.
Such a bundle always exists locally in the sense that each point $t\in T$
has a neighborhood $U$ such that the ideal $\gotoco(U,\A_1\restr U)$ in $A_1$
is strongly Morita equivalent to the ideal $\gotoco(U,\A_2\restr U)$ in $A_2$
via an equivalence that respects $U$; i.e., there is a bundle $\bX_U$
over $U$ so that $\X(t)$ implements a strong Morita equivalence
between $A_1(t)$ and $A_2(t)$ for all $t\in U$. The problem is to glue the $%
\bX_U$ together. The obstruction to doing this is an element in the
sheaf cohomology group $H^2(T,\mathcal{S})$, where $\mathcal{S}$ is the
sheaf of all continuous $\T $-valued functions on $T$. Indeed, given the
continuous trace $\cs$-algebra $A=\gotoco(T,\A)$ one can find an open
cover $\U=\{U_i\}$ of $T$ such that each ideal $\gotoco(U_i,\A
\restr{U_i})$ is strongly Morita equivalent to 
$C_0(U_i)$ over $U_i$. The bundle $%
\bX_{U_i}$ associated with this equivalence is, in fact, a \emph{%
Hilbert bundle}, i.e., the fibres $\X(t)$ are Hilbert spaces, and
one has that $\A\restr{U_i}$ is isomorphic to $\KK(\bX_{U_i})$. 
On overlaps, $U_i\cap U_j$, the bundles ${\KK}({\bX}_{U_i})\restr{U_i\cap
U_j}$ and ${\KK}({\bX}_{Uj})\restr{U_i\cap U_j}$, being strongly Morita
equivalent (and separable), are stably isomorphic. 
Assuming, without loss
of generality, then, that the bundles are stable, we may take the fibres 
${\X}(t)$ to be infinite dimensional. 
We may then infer that there are
unitary valued functions $u_{ij}$ on 
$U_i\cap U_j$ such that $u_{ij}(t):{\X}(t)\rightarrow {\X}(t)$ and
such that $\Ad(u_{ij})$ implements a bundle 
isomorphism from ${\KK}({\bX}_{U_i})\restr{U_i\cap U_j}$ to ${\KK}(%
{\bX}_{U_j})\restr {U_i\cap U_j}$. 
If one compares $u_{ik}$ with $u_{ij}u_{jk}$
on a triple overlap, $U_i\cap U_j\cap U_k$, one sees that there is a
unimodular function $\delta _{ijk}$ such $u_{ik}(t)=\delta
_{ijk}(t)u_{ij}(t)u_{jk}(t)$ on $U_i\cap U_j\cap U_k$.
Straightforward calculations then
reveal that these functions constitute a $2$-cocycle representing an element 
$\delta =\delta (A)$ in $H^2(T,\sT)$. The element $\delta (A)$ is
called the \emph{Dixmier-Douady class} of $A$ and, phrased in modern
terms,
 one of the key results
in \cite{dd} is that two continuous trace $\cs$-algebras with spectrum $T
$ are strongly Morita equivalent over $T$ if and only if their
Dixmier-Douady classes coincide.

Given continuous trace $\cs$-algebras $A_1$ and $A_2$ with spectra $T$,
one can form $A_1\otimes _{C_0(T)}A_2$. In terms of bundles, the bundle for
this tensor product is pointwise tensor product of the bundles $\A%
_1 $ and $\A_2$ representing $A_1$ and $A_2$, respectively. Notice
that if, for $i=1,2$, $A_i^{\prime }$ is strongly Morita equivalent to $A_i$
over $T$, then $A_1\otimes _{C_0(T)}A_2$ is also strongly Morita equivalent
to $A_1^{\prime }\otimes _{C_0(T)}A_2^{\prime }$ over $T$. Thus if, as
above, we let $\br(T)$ be the collection of equivalence classes of continuous
trace $\cs$-algebras with spectrum $T$, where the equivalence relation is
``strong Morita equivalence over $T$'', then $\br(T)$ becomes a semigroup
under the operation of tensoring over $C_0(T)$. In fact, $\br(T)$ is a group,
the Brauer group of $T$: the identity is represented by $C_0(T)$ and the
inverse of $[A]$ is represented by the conjugate of $A$, $\overline{A}$.
That is, $\overline{A}$ is the same set as $A$, with the same product and
addition, but if $\flat :A\rightarrow \overline{A}$ is the identity map,
then $\lambda \cdot \flat (a)=\flat (\overline{\lambda }\cdot a)$, $\lambda
\in \C$, $a\in A$. Alternatively, one may think of $\overline{A}$ as the
opposite algebra of $A$. The main result of \cite{dixfields}, then, is that the
map $[A]\rightarrow \delta (A)$ is an \emph{group isomorphism} of $\br(T)$ with $%
H^2(T,\sT)$.

Now the surjectivity of $\delta $ was proved in \cite{dd} using the fact that
the unitary group of an infinite dimensional Hilbert space is contractible.
Subsequently, in independent investigations, the third author \cite{renault3}
and Raeburn and Taylor \cite{rt}
showed that $\delta $ is onto using groupoids. The idea behind both
approaches, and one that is important for this paper, may be traced back to
Mackey's pioneering article 
\cite[p.~1190]{m6}. In this paper,
Mackey showed how to describe the \emph{topological} cohomology of a space
in terms of the \emph{groupoid} cohomology of certain groupoids constructed
from covers of the space. It is important for our purposes to have some
detail at our disposal about how this is done.

So fix the locally compact Hausdorff space $T$ and a locally finite cover $%
\U=\{U_i\}_{i\in I}$ of $T$. Let $X=\coprod_{i\in I}U_i=\{(i,t):t\in
U_i\}$ and define $\psi :X\rightarrow T$ by $\psi (i,t)=t$. Then, of course, 
$\psi $ is a local homeomorphism and the set $G^\psi :=\coprod_{i,j}U_i\cap
U_j=\{(i,t,j):t\in U_i\cap U_j\}=X*_\psi X$ has a naturally defined
structure of a locally compact groupoid that is $r$-discrete in the sense of
\cite{renault}\footnote{Since we always assume the the range map is
open, a groupoid is $r$-discrete here exactly when the range map is a
local homeomorphism \cite[Proposition~I.2.8]{renault}}. 
The unit space of $G^\psi $ is $X$ and two triples
$(i_1,t_1,j_1)$ 
and $(i_2,t_2,j_2)$ are composable (or multipliable) if and only if $j_1=i_2$
and $t_1=t_2$, in which case $(i_1,t_1,j_1)\cdot
(j_1,t_1,j_2)=(i_1,t_1,j_2)$,
and $(i_1,t_1,j_1)^{-1}=(j_1,t_1,i_1)$. If $\nu =\{\nu _{ijk}\}$ is a
continuous $2$-cocycle for $\U$ with values in $\T $, then $\nu $
determines a continuous groupoid $2$-cocycle $\tilde{\nu}$ on $G^\psi $ with
values in $\T $ by the formula 
\begin{equation*}
\tilde{\nu}((i,t,j),(j,t,k))=\nu _{ijk}(t).
\end{equation*}
And conversely, a continuous $\T $-valued, groupoid $2$-cocycle on $%
G^\psi $ determines a $2$-cocycle on $\U$. Moreover, two cocycles, $%
\nu ^1$ and $\nu ^2$, are topologically cohomologous by a coboundary over $%
\U$ if and only if the associated groupoid cocycles, $\tilde{\nu}^1$
and $\tilde{\nu}^2$, are cohomologous in the groupoid sense. Thus we find
that the groupoid cohomology group $H^2(G,\T)$ is isomorphic to the
sheaf cohomology group $H^2(\U,\sT)$ of the cover~$\U
$.

Just as with groups, we prefer to think of extensions instead of groupoid
cocycles. The groupoid $2$-cocycle $\tilde{\nu}$ gives an extension of $%
G^\psi $ by the $(G^\psi )^{(0)}\times \T $, denoted $E_\nu $.
Set-theoretically, $E_\nu =G^\psi \times \T $, and the product is given
by the formula 
\begin{equation*}
((i,t,j),z)\cdot ((j,t,k),w)=((i,t,k),\nu _{ijk}(t)zw).
\end{equation*}
The topology on $E_\nu $ is the product topology because $\nu $ is
continuous. The extension $E_\nu $ is called a \emph{twist over $G^\psi $}
in \cite{alex}. In general, a twist $E$ over a locally compact groupoid 
$G$ is simply an extension $E$ of $G$ by the groupoid $G^{(0)}\times \T $%
. Such an $E$ has the natural structure of a principal $\T $-bundle over 
$G$ and two twists $E_1$ and $E_2$ are called equivalent if there is a
groupoid isomorphism $\phi :E_1\rightarrow E_2$ that is also a bundle map.
With respect to the operation of forming fibred products, or Baer sums, the
collection of twists over a locally compact groupoid $G$ modulo twist
equivalence is a group, denoted $\Tw(G)$. In the case that the groupoid is $%
G^\psi $, an equivalence $\phi :E_\nu \rightarrow E_\mu $ between two twists
gives rise to continuous functions $\alpha _{ij}:U_i\cap U_j\rightarrow 
\T$ such that 
\begin{equation*}
\phi ((i,t,j),z)=((i,t,j),\alpha _{ij}(t)z),
\end{equation*}
because $\phi $ is a bundle map, and the $\alpha _{ij}$'s further satisfy
the equation 
\begin{equation*}
\mu _{ijk}(t)\alpha _{ij}(t)\alpha _{jk}(t)\overline{\alpha _{ik}(t)}=\nu
_{ijk}(t)
\end{equation*}
because $\phi $ is a homomorphism. Thus the cocycles $\nu $ and $\mu $ are
cohomologous over $\U$. Hence, we conclude that there are natural
isomorphisms between the three groups $H^2(\U,\sT%
)$, $H^2(G^\psi ,\T )$, and $\Tw(G^\psi )$.

Returning to the groupoid approach to showing that the Dixmier-Douady map $%
\delta $ is surjective, discovered by the third author \cite{renault3}
and 
Raeburn and Taylor \cite{rt},
let $[\nu ]$ be a given element of $H^2(T,\mathcal{S})$ and choose a cover $%
\mathbf{U}$ on which $[\nu ]$ is represented by $\nu $. If $E_\nu $ is the
associated twist, then the restricted groupoid $\cs$-algebra $\cs(G^\psi
;E_\nu )$ is continuous trace with spectrum $T$ and $\delta (\cs(G^\psi
;E_\nu ))=[\nu ]$.

The relation between topological cohomology and groupoid cohomology is
further strengthened through the notion of \emph{equivalence of groupoids}
\cite{mrw} and leads to our definition of $\extgt$ for an arbitrary
locally compact groupoid $G$. Indeed, if the space $T$ is viewed as a
groupoid, the so-called \emph{cotrivial} groupoid with unit space $T$, then $%
T$ is equivalent, in the sense of \cite{mrw} to $G^\psi $. Moreover, if $%
\boldsymbol{V}=\{V_j\}_{j\in J}$ is another cover of $T$, if
$Y:=\coprod_{j\in J}V_j 
$, with local homeomorphism $\varphi :Y\rightarrow T$, given by $\varphi
(j,t)=t$, and if $G^\varphi =Y*_\varphi Y$ is the associated groupoid, then $%
G^\psi $ and $G^\varphi $ are equivalent. In fact, $X*Y=\{(x,y)\in X\times
Y:\psi (x)=\varphi (y)\}$ is a $G^\psi ,G^\varphi $-equivalence in the sense
of \cite{mrw}. Notice, however, that $X*Y$ comes from the cover $\U\cap 
\boldsymbol{V}=\{U_i\cap V_j\}_{(i,j)\in I\times J}$ that \emph{refines} $\U%
$ and $\boldsymbol{V}$. Furthermore, if $G^{\psi ,\varphi }=(X*Y)*_{(\psi
,\varphi )}(X*Y)$ is the associated groupoid, then $G^{\psi ,\varphi }$ is
equivalent to both of $G^\psi $ and $G^\varphi $, \emph{and} each of these
groupoids is a homomorphic image of $G^{\psi ,\varphi }$. In this way, we
find that the notion of `refinement' for covers is captured in the notion of
`equivalence' for groupoids. The collection of groupoids equivalent to $T$
becomes directed through the process of taking homomorphisms. We thus arrive
at the isomorphism 
\begin{equation*}
H^2(T,\sT)\cong \underset{\longrightarrow }\lim\; H^2(\U,%
\sT)\cong\underset{\longrightarrow }\lim \Tw(G^\psi )
\end{equation*}
that relates $H^2(T,\sT)$ with groupoid cohomology and suggests how
to proceed in general (see Sections \ref{sec-eh}~and \ref{sec-four}). 

For a fixed locally compact groupoid $G$ with Haar system, we shall consider
the collection $\pg $ consisting of right principal $G$ spaces
(with equivariant $s$-systems in the sense of \cite{renault2}). For each $X\in 
\pg $, we let $G^X$ denote the imprimitivity groupoid associated
with $X$. This is a canonical groupoid with Haar system equivalent to $G$
(see \cite{mrw}). We form the twist group $\Tw(G^X)$ defined above. Owing to the
fact that $X$ can carry non-trivial first cohomology, we actually take a
quotient of $\Tw(G^X)$ reflecting this, getting a group that we call $%
\E(G^X)$. As we shall see, $\pg $ forms a directed system
and the family $\{\E(G^X)\}_{X\in \pg }$ is inductive. We
define, then, $\extgt$ as the limit 
\begin{equation*}
\underset{X\in \pg }{\underset{\longrightarrow }{\lim }}\E%
(G^X).
\end{equation*}
This is our way to generalize $2$-cohomology of a space to groupoids. It
looks like one could develop a whole cohomology theory along these lines,
but at the moment, the technical difficulties seem formidable. In any case,
this $\extgt$ suffices for our purposes here.

As we indicated above, $\br(T)$ may be viewed in terms of strong
Morita equivalence classes of elementary $\cs$-bundles over $T$, each
satisfying Fell's condition. To generalize this notion to groupoids and to
return to the Bargmann-Mackey-Wigner analysis, we proceed as follows. Fix a
locally compact groupoid $G$ with Haar system. We consider pairs $(\A%
,\alpha )$ consisting of elementary $\cs$-bundles $\A$ over $%
G^{(0)}$ satisfying Fell's condition and continuous actions $\alpha $ of $G$
on $\A$. Of course, if $G$ is a group, the bundles reduce to a
single elementary $\cs$-algebra, i.e., the algebra of compact operators on
a Hilbert space, and we are looking at its automorphic representations. We
shall say that two pairs $(\A,\alpha )$ and $(\B,\beta )$
are \emph{Morita equivalent} if there is an $\A\sme\B$%
-imprimitivity bimodule bundle $\bX$ over $G^{(0)}$ that admits a
continuous action $V$ of $G$ by isometric $\C $-linear isomorphisms such
that $V$ implements the actions $\alpha $ and $\beta $ in a natural way (see
Definition~\ref{def-3.1}). 
In the case when $G$ is a group and the bundles reduce to $\K
$, this notion reduces to \emph{exterior equivalence} described above. (See
Section~\ref{sec-e(x)}, below.) The Morita equivalence class of $(\A,\alpha )$
will be denoted $[\A,\alpha ]$ and the collection of all such
equivalence classes will be denoted by $\brg$. Under the operation of
tensoring, as we shall see, $\brg$ becomes a group. The identity element is
represented by the class of $(G^{(0)}\times \C ,\tau )$, where $\tau
_\gamma (s(\gamma ),c)=(r(\gamma ),c)$, and the inverse of $[\A%
,\alpha ]$ is represented by $[\overline{\A},\overline{\alpha }]$,
where $\overline{\A}$ is the bundle of conjugate or opposite $\cs$%
-algebras of $\A$, and where $\overline{\alpha }$ is $\alpha $. The
details of the proof of this theorem, \thmref{thm-grp}, take some time to
develop, Sections \ref{sec-one}~and \ref{sec-bg}, 
but the proof is conceptually straightforward. It
follows the lines of Theorem 3.6 in \cite{ckrw}.

It should be emphasized that our $\brg$ coincides with $\br(T)$ when $G$
reduces to a space $T$ and that when $G$ is the transformation group
groupoid, $X\times H$, determined by a locally compact group $H$ acting on
locally compact space $X$, then our $\brg$ coincides with the \emph{%
equivariant Brauer group} $\br_H(X)$ analyzed in \cite{ckrw}. Of course, in the
very special case when $X$ reduces to a point, so that $G$ may be viewed as
the group $H$, we find that 
$\brg=\br_H(\set{pt})\cong \mathcal{R}(G)/\mathcal{S}
$ described above.

Our ultimate goal is to prove \thmref{thm-main3}, which asserts that $\brg\cong
\extgt$. As we have said, this result contains, simultaneously, the
isomorphism of $\brg\;(=\mathcal{R}(G)/\sT)$ with the Borel
cohomology group $H^2(G,\T )$, when $G$ is a locally compact group, and
the isomorphism $\br(T)\cong H^2(T,\sT)$ when $T$ is a space.

When the groupoid $G$ is the transformation group groupoid, $X\times H$,
Crocker and Raeburn, together with the first and fourth authors, described $%
\brg\;(=\br_H(X))$ in terms of certain Moore cohomology groups (\cite{moore3}
and \cite{moore4}) arising from various actions of $H$ on the objects that
are involved in the analysis of $\cs$-dynamical systems coming from $H$ and
continuous trace $\cs$-algebras with spectrum $X$ \cite{ckrw}. We shall not
enter into the details here, except to say that the Moore cohomology groups
are measure theoretic objects, like Mackey's $H^2(G,\T )$. Our approach
gives a purely \emph{topological} description $\brg$. We concentrate on
extensions, as emphasized above, and not on actions of $H$ on Polish
groups, as in \cite{ckrw}. However, it remains to be seen what the precise
relation is between our description of $\brg$ and that in \cite{ckrw}. In
particular, it would be interesting to understand the ``filtration''
for $\brg$, discovered there, in terms of groupoid theoretic constructs.
Also, as noted above, our analysis suggests that a cohomology
theory based on $\mathcal{P}(G)$ might exist that 
generalizes the measure theoretic cohomology of Moore.

The first
author \cite{alex3} adapted the equivariant sheaf cohomology theory of
Grothendieck \cite{groth} to cover an $r$-discrete groupoid $G$
acting on sheaves defined on the unit space $G^{(0)}$, and he showed
that  $\brg$ is isomorphic to an equivariant
sheaf cohomology group $H^2(G,\sT)$. Here $\sT$ is the sheaf
of germs of continuous $\T $-valued functions on the unit space, $G^{(0)}
$, of $G$, with $G$ acting on $\sT$ in a natural way. The
equivariant sheaf cohomology was invented in part to obtain a deeper
understanding of the group $\Tw(G)$, which is an extension group. Our work
gives a meaning to $H^1(G,\sT)$ and $H^2(G,\sT)$, at least,
when $G$ is not $r$-discrete, but at the present time, we are not able to go
beyond this. We believe that rationalizing and unifying the various
cohomology theories that we have been discussing is one of the most
important challenges for the theory in the future.

The paper is organized as follows. In Section~\ref{sec-one}, 
we present preliminaries on
Banach bundles and related material. We define and discuss Morita
equivalence of $\cs$-bundles, imprimitivity bimodule bundles, and actions
of groupoids on bundles. We analyze in detail $\sbr(G)$, the collection
of all $\cs$-$G$-bundles $(\A,\alpha )$, where $G$ is a groupoid
and $\A$ is an elementary $\cs$-bundle over $G^{(0)}$ that is
acted upon by $G$ via $\alpha $. Most important, we present the fundamental
construction that assigns to each principal $G$-space $X$ and each $\cs$-$G
$-bundle, $(\A,\alpha )\in \sbr(G)$, an \emph{induced} $\cs$-%
$G^X$-bundle $(\A^X,\alpha ^X)\in \sbr(G^X)$, where $G^X$ is
the imprimitivity groupoid of $G$ determined by $X$. 

In the next section we define what it means for two systems in $\sbr(G)$
to be Morita equivalent and we define  $\brg$ to be $\sbr(G)$ modulo
this equivalence. We then show that $\brg$ is a group. As we have
indicated, the details follow the lines found in \cite{ckrw}. 

The analysis of $\brg$ really is about three isomorphism theorems. The
first is presented in Section~\ref{sec-fiso}, 
where we present an explicit isomorphism $%
\phi ^X$ between $\brg$ and $\br(H)$ where $X$ is an equivalence between $G$
and $H$. This isomorphism theorem yields as an immediate corollary the main
theorem in \cite{krw}.

Section~\ref{sec-haar} is 
devoted to the relation between Haar systems on equivalent
groupoids, on the one hand, and the notion of equivariant $s$-systems in the
sense of \cite{renault4}, on the other. 
It turns out that the two notions are essentially the same and
this fact plays an important r\^ole in our analysis.

Section~\ref{sec-graph} 
discusses the relation between the notions of groupoid equivalence
and groupoid homomorphisms. In a sense it gives us a start to the
generalization, to arbitrary locally compact groupoids, of the order
structure on $\pg $ that is embodied in the concept of refinement
for covers of a topological space. Most important for our purposes is the
result, \corref{cor-trivialdd}, 
that allows us to replace a system $(\A%
,\alpha )$ in $\sbr(G)$ by the system $(\A^X,\alpha ^X)$ in $%
\sbr(G^X)$ for a suitable principal $G$-space $X$ such that the
Dixmier-Douady class $\delta (\A^X)$ vanishes. 

In Section~\ref{sec-eh}, 
we analyze twists over $G$, but we use ideas from \cite{renault3} and
\cite{gcddi} 
to replace the equivariant sheaf cohomology of \cite{alex3} that was used
for this purpose when $G$ is $r$-discrete. As we have implied, the sheaf
cohomology theory does not seem to have a serviceable analogue in our more
general setting. However, we still are able to show that the group of twists
over $G$, $\Tw(G)$, contains a subgroup of twists that is ``parameterized''
by principal circle bundles on $G^{(0)}$, i.e., by a certain $H^1$ group.
The quotient of $\Tw(G)$ by this subgroup is denoted $\E(G)$ and is
called the generalized twist group. 

With $\E(G)$ in hand, we are able to follow the first author's lead
in \cite{alex3}, where $\brg$ is analyzed in terms of the subgroup $\brz(G)$
consisting of those equivalence classes of systems $(\A,\alpha )$
with trivial Dixmier-Douady invariant. For an equivalence class of $(%
\A,\alpha )$ in $\brz(G)$ one may define a natural twist $E(\alpha )
$ over $G$. The definition is essentially Mackey's definition of $E(\alpha )$
in the context of automorphic representations of locally compact groups
described above. Section~\ref{sec-e(x)} 
is devoted to showing that every twist over $G$
is an $E(\alpha )$ for a suitable $\alpha $,
\propref{prop-surjectivity}, 
and to showing
that $\brz(G)$ is naturally isomorphic to $\E(G)$, Corollaries 
\ref{cor-injhm}~and \ref{cor-singlecover}. 
This last result is our second isomorphism theorem.

Section~\ref{sec-four} is 
devoted to the definition of $\extgt$ as an inductive
limit of the generalized twist groups $\E(G^X)$, where $X$ runs
over the collection of right principal $G$ spaces with equivariant $s$%
-systems, $\pg$. The problem here is to show that $\pg$
is directed. The solution is to rephrase the notion of ``refinement'' for
covers in terms of equivalence of groupoids. The details require technology
developed in Sections~\ref{sec-graph}~and \ref{sec-eh}.

Finally, our third, and main, isomorphism theorem, \thmref{thm-main3}, that
asserts that $\brg$ is isomorphic to $\extgt$, is proved by using
our first and second isomorphism theorems to realize embeddings of the $%
\E(G^X)$, $X\in \pg $, in a fashion that is compatible
with the order on $\pg $. That is, we show that $\brg$ satisfies
the universal properties of the inductive limit defining $\extgt$.
 
The last section is devoted to discussions of a number of examples that help
to illustrate our theory and to establishing a connection between it and the
third author's work on dual groupoids.

\else

\fi
\ifdraft
%
%

\section{Preliminaries}\label{sec-one}

We will need to make considerable use of the notion of a \emph{Banach
bundle}.  An excellent summary of the definitions and results needed
can be found in \S\S13--14 of Chapter~II of \cite{fell-doran} (see
also \cite[\S1]{fell4}).  
We record some of the basic
definitions and results here for the sake of completeness.

\begin{definition}
\label{def-bb}
Suppose that $T$ is a locally compact space.
A \emph{Banach bundle over $T$} is a topological space $\A$ together
with a continuous, open surjection $p=\pa:\A\to T$ and Banach space
structures on each fibre $p^{-1}\(\set t\)$ ($t\in T$) satisfying
the following axioms.
\begin{enumerate}
\item
The map $a\mapsto\|a\|$ is continuous from $\A$ to $\R^+$.
\item
The map $(a,b)\mapsto a+b$ is continuous from 
\[
\A*\A:=\set{(a,b)\in\A\times\A: p(a)=p(b)}
\]
to $\A$.
\item
For each $\lambda\in\C$, the map $b\mapsto\lambda b$ is continuous
from $\A$ to $\A$.
\item
If $\set{a_i}$ is a net in $\A$ such that $p(a_i)\to t$ in $T$ and
such that $\|a_i\|\to 0$ in $\R$, then $a_i\to 0_t$ in $\A$.  
(Of course, $0_t$ denotes the zero element in the Banach space
$p^{-1}\(\set 
t\)$.)
\end{enumerate}
A \emph{\cs-bundle over $T$} is a Banach bundle $p:\A\to T$ such that
each fibre is a \cs-algebra satisfying, in addition to axioms
\partref1--\partref4 above, the following axioms.
\begin{enumerate}
\setcounter{enumi}{4}\item
The map $(a,b)\mapsto ab$ is continuous from $\A*\A$ to $\A$.
\item
The map $a\mapsto a^*$ is continuous from $\A$ to $\A$.
\end{enumerate}
A continuous function $f:T\to \A$ such that $p\circ f=\id_T$ is
called a \emph{section} of $\A$.  The collection of sections $f$ for
which $t\mapsto \|f(t)\|$ is in $C_0(T)$ is denoted by
$A:=\gotoco(T,{\A})$.
\end{definition}

\begin{remark}
If $\A$ is a Banach bundle (resp.{} \cs-bundle) over $T$, then
$A=\gsuboT(\A)$ is a Banach space (resp.{} \cs-algebra) with respect
to the obvious pointwise operations and norm $\|f\|:=\sup_{t\in
T}\|f(t)\|$.
We will use the notation $A(t)$ for the fibre $p^{-1}\(\set t\)$.
Furthermore, $\gotoco(T,\A)$ is a central $C_b(T)$-bimodule.
\end{remark}

\begin{definition}
A \cs-bundle is called \emph{elementary} if every fibre is an
elementary \cs-algebra; that is, every fibre is isomorphic to the
compact operators on some complex Hilbert space.
\end{definition}

\begin{definition}
If $p:\A\to T$ is a Banach bundle and if $q:X\to T$ is continuous,
then the \emph{pull-back of $\A$ along $q$} is the Banach bundle
$p':q^*\A\to X$, where  
\[
q^*\A:=\set{(x,a)\in X\times\A: q(x)=p(a)},
\]
and $p'(x,a):=x$.
\end{definition}

\begin{remark}
Fell refers to $q^*\A$ as the \emph{retraction} of $\A$ by $q$.
Note that the fibre of $q^*\A$ over $x$ can be identified with the
fibre of $\A$ over $q(x)$.  We will often make this identification
without any fanfare.
\end{remark}

If $p:\A\to T$ is a Banach -bundle and $A=\gsuboT(\A)$, then will
often write $q^*(A)$ for $\gsuboT(q^*\A)$.  If $\A$ is a \cs-bundle,
then this is the same notation
employed for the \cs-algebraic pull back
\[
C_0(X)\tensor_{C(T)} A.
\]
This should not be a source of confusion since the two are isomorphic
by \cite[Proposition~1.3]{rw}.

Finally, if each fibre of $\A$ is a Hilbert space, then we
refer to
$\A$ as a \emph{Hilbert bundle}.  A straightforward
polarization argument implies that $(a,b)\mapsto (a\mid b)$ is
continuous from
$\A*\A$ to $\C$.  Therefore $A=\gsuboT(\A)$ is a Hilbert
$C_0(T)$-module \cite{lance}.  This is in contrast to a \emph{Borel
Hilbert bundle}, where the (Borel) sections form a Hilbert
space.

\begin{definition}
A Banach bundle $p:\A\to T$ has \emph{enough cross sections} if given
any $a\in \A$, there is a cross section $f\in \gsuboT(\A)$ such
that $f\(p(a)\)=a$.
\end{definition}

\begin{remark}
Since we are assuming that our base spaces are locally compact, it
follows from a result of Douady and dal Soglio-H\'erault
\cite{douady-herault} that \emph{all} our bundles have enough cross
sections.  A proof of their result appears in Appendix~C of
\cite{fell-doran}.
\end{remark}

\begin{remark}
\label{rem-lee}
It follows from \cite[Theorem~10.5.4]{dix} and Fell's
characterization of Banach bundles (e.g., \cite[Theorem
II.13.8]{fell-doran}) that a CCR \cs-algebra has Hausdorff
spectrum $T$ if and only if $A$ is isomorphic to the section
algebra of an elementary \cs-bundle $\A$ over $T$.
\end{remark}

\begin{remark}
Axiom \partref4 of Definition~\ref{def-bb} ties the topology on a
Banach bundle to its sections.  Specifically, one can use
\cite[Proposition~II.13.12]{fell-doran} to characterize the topology as
follows. 
If $p:\A\to T$ is a Banach bundle, and if $A=\gsuboT(\A)$, then a
net $\set{a_i}$ converges to $a$ in $\A$ if and only if
(a)~$\set{p(a_i)} \to p(a)$, and~(b) $\bigl\|a_i-f\(p(a_i)\)\bigr\|
\to \bigl\|a- f\(p(a)\)\bigr\|$ for all $f\in A$.  If fact, for the
``if direction'', it suffices to take $f\in A$ with $f\(p(a)\)=a$.
\end{remark}

\begin{remark}
\label{rem-**}
Suppose that $\A$ and $\B$ are Banach bundles over a locally compact
space $T$, and that $\phi:\gotoco(T,\A)\to\gotoco(T,\B)$ is an isometric
$C_0(T)$-isomorphism; i.e., $\phi$ respects the $C_0(T)$-module structures. 
Then we obtain a well defined linear map
$\Phi:\A\to \B$ by $\Phi\(f(x)\)=\phi(f)(x)$.  This is a bundle
isomorphism by \cite[Proposition~II.13.16]{fell-doran}.
\end{remark}

Of course, the most studied sorts of bundles are those which are
locally trivial.  However, the bundles which arise from \cs-algebras
(see Remark~\ref{rem-lee}) are almost never naturally locally
trivial, as the fibres are formally distinct spaces.  
Instead the proper notion is the following.

\begin{definition}
\label{def-2.11}
An elementary \cs-bundle over $T$ satisfies \emph{Fell's condition} if
each $s\in T$ has a \nbhd{} $U$ such that there is a section $f$ with
$f(t)$ a rank-one projection for all $t\in U$.
\end{definition}

\begin{remark}
\label{rem-ct}
An elementary \cs-bundle satisfies Fell's condition if and only if its
section algebra is a continuous-trace \cs-algebra
\cite[Proposition~10.5.8]{dix}. 
\end{remark}

We will need to use the notion of a groupoid action on a space.
Let $G$ be a groupoid and let $X$ be
a topological space together with a continuous open map $s$
from $X$ onto $\go$. 
Form $X*G=\set{(x,\gamma):s(x)=r(\gamma)}$. (We have chosen not to
write $s_X$, $s_G$, etc., 
as these notations eventually become too distracting.)
To say that $G$ acts (continuously) on the right of $X$
means that there is a
continuous map from $X*G$ to $X$, with the image of $(x,\gamma)$
denoted $x\cdot\gamma$, such that the
following hold:
\begin{enumerate}
\item
$s(x\cdot\gamma)=s(\gamma)$,
\item
$x\cdot(\alpha\beta)=(x\cdot\alpha)\cdot\beta$, for all $(\alpha,
\beta)\in G^{(2)}$, and
\item
$x\cdot s(x)=x$ for all $x\in X$.
\end{enumerate}
We think of $ s$ as being a ``generalized source'' map for the
action.  Likewise a (continuous) left action of $G$ on $X$ is
determined by a continuous open surjection $r:X\to\go$, and a
continuous map from $G*X$ to $X$ satisfying the appropriate
analogues of \partref1, \partref2,~and \partref3.
In a like fashion, we think of $r$ as a ``generalized range''
map.  When we speak of left or right actions in the sequel, the
maps $ s$ and $r$ will be implicitly understood and will
be referred to without additional comment.  Most of our
definitions will be made in terms of right actions, but the left
handed versions can be formulated with no difficulty.

As with group actions, if $G$ acts on $X$ (on the right, say),
and if $x\in X$, the {\it orbit\/} of $x$ is simply $\set {
x\cdot\gamma:\gamma\in G^{s(x)}}$.  The orbits partition $X$
and we write $X/G$ for the quotient space with the quotient
topology.

Consider a (right) action of $G$ on $X$ and let
$\Psi:X*G\to X\times X$ be defined by the formula
$\Psi(x,\gamma)=(x,x\cdot\gamma)$.  Then we call the action {\it
free\/} if $\Psi$ is one-to-one.  Alternatively, the action is
free precisely when the equation $x\cdot\gamma=x$ implies that
$\gamma=s(x)$.  The action is called {\it proper\/} provided
that $\Psi$ is a proper map.

Now suppose that $G$ is a locally compact groupoid  
and that $p:\A\to\go$
is a Banach bundle (resp.{}
\cs-bundle).  A collection $(\A,\alpha)$ of isometric isomorphisms
(resp.{} $*$-isomorphisms) $\alpha_\gamma:A\(s(\gamma)\)\to
A\(r(\gamma)\)$ is called a \emph{$G$-action by isomorphisms} or
simply a $G$-action on $\A$ if $\gamma\cdot a:=\alpha_\gamma(a)$ makes
$(\A,p\)$ into a (topological) left $G$-space as in
\cite[\S2]{gcddi}.  In this case, we say that $\A$ is a Banach
$G$-bundle (resp.{} \cs-G-bundle). 

If $\alpha$ is a $G$-action on $\A$, then $\alpha$ induces a
$C_0(\go)$-linear isomorphism $\alpha:s_G^*(A)\to r_G^*(A)$ defined by
$\alpha(f)(\gamma):=\alpha_\gamma\(f(\gamma)\)$.  It is not hard to
see that the converse holds.

\begin{lem}
\label{lem-1.12}
Let $G$ be a locally compact groupoid.
Suppose that $p:\A\to \go$ is a Banach bundle, and that for each $\g\in
G$, $\ag:A\(s(\g)\)\to A\(r(\g)\)$ is an isometric isomorphism
such that $\alpha_{\g\eta}=\alpha_\g\circ\alpha_\eta$ if
$(\g,\eta)\in G^{(2)}$.   Then
$(\A,\alpha)$ is a $G$-action by isomorphisms if and only if
$\alpha(f)(\g) := \ag\(f(\g)\)$ defines an $C_0(\go)$-isomorphism of
$s_G^*(A)$ 
onto $r_G^*(A)$. 
\end{lem}

\begin{proof} 
We have already observed that if $\alpha$ is a $G$-action, then
$\alpha:s_G^*(A) \to r_G^*(A)$ is a $C_0(\go)$-isomorphism.

To see the converse, the only issue is to observe that the action is
continuous.   Suppose that $(\gamma_i,a_i)\to (\gamma,a)$
in $s_G^*(\A)$.  Then there is a $g\in s_G^*(A)$ such that
$g(\gamma)=a$.  Since $\alpha(g)\in r_G^*(A)$, $\alpha(g)(\g_i) \to
\alpha(g)(\g) = \g\cdot a$.  On the other hand, 
\begin{equation*}
\|\alpha(g)(\g_i) - \g_i\cdot a\|  = \| \alpha_{\g_i}\( g(\g_i) - a_i
\)\| 
= \| g(\g_i) - a_i\|,
\end{equation*}
which converges to zero.  It now follows from
Definition~\ref{def-bb}\partref4 that $\g_i\cdot a_i\to \g\cdot a$ as
required.
\end{proof}

\begin{definition}
\label{def-2.14}
Let $G$ be a \sclcg.
  We denote by $\sbrg$ the collections of
\cs-$G$-bundles
$(\A,\alpha)$ where $\A$ is an elementary \cs-bundle with
separable fibres and which satisfies Fell's condition.\footnote{Equivalently,
$A=\gsuboT(\A)$ is a separable continuous-trace 
\cs-algebra.}
\end{definition}

Now we turn to a construction which will be crucial in the sequel.
Let $X$ be a $(H,G)$-equivalence\footnote{For the construction, it
will suffice for $X$ to be a principal right $G$-space and a left $H$
space such that $(h\cdot x)\cdot\g=h\cdot(x\cdot\g)$ and such that
$r_X:X\to\ho$ induces a homeomorphism of $X/G$ with $\ho$.}.  Then
$s_X^*\A = \set{(x,a)\in X\times\A:s_X(x)=p(a)}$ is a (not necessarily
locally compact) principal right $G$-space:
\[
(x,a)\cdot\g:=\(x\cdot\g,\ag^{-1}(a)\).
\]
Define $\ax$ to be the quotient space $s_X^*\A/G$.  We will employ the
notation $[x,a]$ to denote the orbit of $(x,a)$ in $\ax$.  The map
$(x,a)\mapsto r_X(x)$ defines a continuous surjection $\px:\ax\to
\ho$.  To see that $\px$ is open\footnote{One useful criterion for a
map to be open is that convergent nets in the image have a subnet
which can be lifted to the domain.  An explicit statement and proof
can be found, for example, in \cite[Proposition~II.13.2]{fell-doran}},
suppose that $\rx(x_i)\to \rx(x)$ in $\ho$.  Passing to a subnet and
relabeling, we can assume that there are $\g_i$ in $G$ such that
$x_i\cdot\g_i\to x$ in $X$.  Since $p:\A\to\go$ is open, we may as
well assume that there are $a_i\in\A$ such that $a_i\to a$ in $\A$ and
$p(a_i)=s(\g_i)$.  Then $(x_i\cdot\g_i,a_i)\to (x,a)$ in $s_X^*\A$.
The openness of $\px$ follows from this.
Finally, the map $a\mapsto [x,a]$ defines an isomorphism of
$A\(\sx(x)\)$ onto $A^X\(\rx(x)\):= (\px)^{-1}\(\rx(x)\)$.  (Note that
the operations and norm on  $(\px)^{-1}\(\rx(x)\)$ are independent of
our choice of representative for $\rx(x)$.)

\begin{prop}
\label{prop-fundcon}
Let $G$ and $H$ be locally compact groupoids.
Suppose that $X$ is a $(H,G)$-equivalence and that $\A$ is a Banach
$G$-bundle.
Then $\px:\ax\to\ho$ is a Banach bundle.
Futhermore, $\a^X_h[x,a]:= [h\cdot x,a]$
defines a left $H$-action on $\ax$ making $(\ax,\a^X)$ an Banach
$H$-bundle called the \emph{the induction of $A$ from $G$ to
$H$ via $X$}. In particular, if
$(\A,\a)\in\sbrg$, then $(\ax,\a^X)\in\sbr(H)$.
\end{prop}

\begin{proof}
For the first assertion we must verify axioms \partref1--\partref4 of
Definition~\ref{def-bb}.  For example, suppose that $[x_i,a_i]
\to [x,a]$ in $\ax$.  Note that $\|(x_i,a_i)\|=\|a_i\|$.
Thus, if \partref1 fails, then we can pass to a subnet and relabel so
that $\bigl|\|a_i\|-\|a\|\bigr|\ge\epsilon>0$ for all $i$.  We can
assume that there are $\g_i\in G$ such that $\(x_i\cdot
\g_i,\a_{\g_i}^{-1}(a_i)\) \to (x,a)$ in $s_X^*\A$.  This leads to a
contradiction since $\|a_i\|=\|\a_{\g_i}^{-1}(a_i)\|$.

To establish \partref2, suppose that $[x_i,a_i]\to [x,a]$ and 
$[y_i,b_i]\to
[y,b]$ in $\ax$ with
$r_X(x_i)=r_X(y_i)$ for all $i$.  We need to show that $[x_i,a_i] +
[y_i,b_i]
\to [x,a]
 + [y,b]$.  But $y_i=x_i\cdot \g_i$ for some $\g_i\in G$ with
$s_G(\g_i) = p(b_i)$.  Thus
\[
[x_i,a_i]+ [y_i,b_i]= [x_i,a_i]+
[x_i\cdot\g_i,b_i]= [x_i,a_i + \a_{\g_i}(b_i)].
\]
It results that we can assume that $y_i=x_i$ in the following.
Therefore if \partref2 fails, we can pass to a subnet and relabel so
that there is a \nbhd{} $U$ of $[x,a+b]$ which contains no
$[x_i,a_i+b_i]$.  But we may as well also assume that there are
$\gi,\eta_i\in G$ such that
$\(x_i\cdot\gi,
\a_{\gi}^{-1}(a_i)\) \to (x,a)$ and
$\(x_i\cdot\eta_i,\a_\gi^{-1}(b_i)\) \to (x,b)$.  Since $X$ is a
principal $G$-space, we may pass to yet another subnet and relabel so
that we can assume $\eta_i^{-1}\gi\to \sx(x)$ in $G$.  Thus
$\(x_i\cdot\gi, \a_\gi^{-1}(b_i)\)\to (x,b)$, and $\(x_i\cdot\gi ,
\a_\gi^{-1}(a_i) + \a_\gi^{-1}(b_i)\)\to (x,a+b)$.  This implies that
$[x_i,a_i+b_i]\to [x,a+b]$, which contradicts our
assumptions.  Axioms \partref3 and \partref4 follow from similar (but
easier) arguments.  

Similarly, if $\A$ is a \cs-bundle, then $\ax$ also satisfies axioms
\partref5 and \partref6; that is, $\ax$ is a \cs-bundle which is
clearly elementary if $\A$ is.

It is straightforward to verify that $\a^X$ is an
action of $H$ on $\ax$ by isomorphisms (or $*$-isomorphisms if $\A$ is
a \cs-bundle).  It only remains to check that $\ax$ satisfies Fell's
condition if $\A$ does.  However, the map $(x,a)\mapsto \(x,(x,a)\cdot
G\)$ 
is a Banach bundle isomorphism\footnote{If $\A$ and $\B$ are bundles
over $T$, then a continuous map $\psi:\A\to\B$ is a \emph{bundle
map} if the diagram
\[
\begin{diagram}
\node{\A}\arrow[2]{e,t}{\psi} \arrow{se,b}{\pa} \node[2]{\B}
\arrow{sw,b}{\pb} \\
\node[2]{T}
\end{diagram}
\]
commutes.  We say $\psi$ is a bundle isomorphism if
$\psi$ is a bijection, and $\psi$ is bicontinuous.  If $\A$ and $\B$ are Banach
bundles, then a bundle isomorphism $\psi$ is a Banach (resp.{} \cs-) bundle
isomorphism if it is also an isomorphism (resp.{} $*$-isomorphism) 
of the fibres.  Note that a continuous bundle map $\psi$ which is an
isomorphism of each fibre is necessarily a Banach (resp.{} \cs-)
bundle isomorphism \cite[Proposition~II.13.17]{fell-doran}.} of
$\sx^*\A$ onto $\rx^*\ax$.  Since $\A$ is assumed to satisfy Fell's
condition, it is clear that $\sx^*\A$ does.  Consequently, so does
$\rx^*\ax$.  Therefore $\rx^*(A^X)=\gotoco(X,\rx^*\ax)$ has continuous
trace (Remark~\ref{rem-ct}).  Since $\rx:X\to\ho$ is continuous, open,
and surjective, \cite[Lemma~1.2]{rr} implies that $A^X$ has continuous
trace.  Therefore $\ax$ satisfies Fell's condition as required.
\end{proof}

If $A$ and $B$ are \cs-algebras and if $\X$ is an
$A\sme B$-\ib, then there is a homeomorphism $h_\X:\hat A\to \hat B$.
If the spectrums of $A$ and $B$ have been identified with $T$, there
is no reason that $h_\X$ has to be the identity map.  When $A$ and $B$
have spectrum (identified with) $T$, then we call $\X$ a $A\smeover
TB$-\ib{} when $h_\X=\id_T$, and we say that $A$ and $B$ are Morita
equivalent over $T$ (cf., \cite[p.~1035]{90c}).  Note that if
$p:\A\to T$ is an elementary \cs-bundle, the spectrum of
$A=\gsuboT(\A)$ is naturally identified with $T$.

\begin{definition}
Two \cs-bundles $\pa:\A\to T$ and $\pb:\B\to T$ with section algebras
$A$ and $B$, respectively, are \emph{Morita
equivalent} if there exists an $A\smeover TB$-\ib{} $\X$.
\end{definition}

It follows from \cite[Theorem~II.13.18]{fell-doran} and
\cite[Corollary~II.14.7]{fell-doran} that $\X=\gsuboT(\bX)$ for a
unique Banach bundle $q:\bX\to T$.  Notice that the fibres $\X(t)$ are
exactly the quotient $A(t)\sme B(t)$-\ib s as described in
\cite[Corollary~3.2]{rieff2}.  In particular if $f\in A$, $g\in B$,
and $\xi,\eta\in\X$, then we have the following:
\begin{equation}
\label{eq-***}
\begin{split}
\lip A<\xi,\eta>(t) &= \blip A(t) <\xi(t),\eta(t)>, \\
\rip B<\xi,\eta>(t) &= \brip B(t)  <\xi(t),\eta(t)>, \\
(f\cdot \xi)(t) &= f(t)\cdot \xi(t)\text{, and} \\
(\xi\cdot g)(t) &= \xi(t)\cdot g(t).
\end{split}
\end{equation}
It follows that there are continuous maps $(x,y)\mapsto \lip\A<x,y>$
from $\bX*\bX:=\set{(x,y)\in\bX\times\bX:q(x)=q(y)}$ to $\A$ and
$(x,y) \mapsto \rip\B<x,y>$ from $\bX*\bX\to \B$ such that, for
example,
$\lip\A<\xi(t),\eta(t)>=\lip A(t) <\xi(t),\eta(t)>$.  Similar, there
are continuous maps $(a,x)\mapsto a\cdot x$ from
$\A*\bX=\set{(a,x)\in\A\times \bX: \pa(a)=q(x)}$ to $\bX$, and
$(x,b)\mapsto x\cdot b$ from $\bX*\B$ to $\bX$.  

\begin{definition}
Suppose that $\pa:\A\to T$ and $\pb:\B\to T$ are \cs-bundles
with section algebras
$A$ and $B$, respectively.  Then a Banach bundle $q:\bX\to T$ is
called an $\A\sme \B$-\ibb{} if each fibre $\X(t)$ is an $A(t)\sme
B(t)$-\ib{} 
such that the natural maps $(a,x)\mapsto a\cdot x$ from $\A*\bX$ to
$\bX$ and $(x,b)\mapsto x\cdot b$ from $\bX*\B$ to $\bX$ are continuous.
\end{definition}

As we will be working mostly with Banach bundles, it will be more
natural here to work with \ibb s rather than \ib s.  Fortunately our
next result will make it easy to go from one to the other.

\begin{prop}
Suppose that $\pa:\A\to T$ and $\pb:\B\to T$ are \cs-bundles
with section algebras
$A$ and $B$, respectively, and that $q:\bX\to T$ is an $\A\sme
\B$-\ibb.  Then $\X:=\gsuboT(\bX)$ is an $A\smeover TB$-\ib{} with
actions and inner products defined as in \eqref{eq-***}.
Conversely, if $\X$ is an $A\smeover TB$-\ib, then there is a unique $\A\sme
\B$-\ibb{} $\bX$ such that $\X\cong \gotoco(T,\bX)$.
\end{prop}
\begin{proof}
The last assertion follows from the preceding discussion.  On the
other hand, if $\X=\gotoco(T,\bX)$, then
it is immediate that $\X$ is an $A\sme B$-bimodule.  Next we observe
that $(x,y)\mapsto \rip\B<x,y>:= \rip B\(q(x)\)<x,y>$ is continuous
from $\bX*\bX$ to $\B$; this follows from a polarization argument and
the fact that $x\mapsto \|x\|=\rip\B<x,x>$ is continuous since $\bX$
is a Banach bundle.  The Cauchy-Schwartz inequality for Hilbert
modules implies that $\|\rip\B<x,y>\|\le\|x\|\|y\|$.  Consequently
\[
\rip B<\xi,\eta>(t):= \brip\B<\xi(t),\eta(t)>
\]
defines a $B$-valued inner product on $\X$.  Then ideal
$\rip\B<\X,\X>$ is dense in $B$ as it is a $C_0(T)$-submodule and
$\rip \B<\X,\X>(t)$ is dense in $B(t)$ for each $t\in T$.  Similar
arguments apply to the $A$-valued inner product.  The rest is
straightforward. 
\end{proof}
\else

\fi
\ifdraft
%
%

\section{The Brauer Group}\label{sec-bg}

The notion of Morita equivalence for \cs-$G$-systems that we are about
to define is crucial for all the work to follow.  It is a natural
generalization to groupoids of the idea of Morita equivalence for
\cs-dynamical systems first investigated in \cite{combes} and
\cite{cmw}.

\begin{definition}
\label{def-3.1}
Two \cs-$G$-systems $(\A,\a)$ and $(\B,\beta)$ are \emph{Morita
equivalent} if there is a $\A\sme\B$-\ibb{} $\bX$ which admits an
action $V$ of $G$ by isomorphisms such that
\begin{gather*}
\blip\A<\vg(x),\vg(y)> = \alpha_\gamma\(\lip \A<x,y>\),\text{ and}\\
\brip\B<\vg(x),\vg(y)> = \beta_\g \(\rip\B<x,y>\).
\end{gather*}
In this case we will write $(\A,\a)\sim_{(\bX,V)}(\B,\beta)$.
\end{definition}

\begin{lem}
Morita equivalence of \cs-$G$-systems is an equivalence relation.
\end{lem}
\begin{proof}
Since $(\A,\a)\sim_{(\A,\a)}(\A,\a)$, ``$\sim$'' is reflexive.  Recall
that if $\X$ is an $A\sme B$-\ib, then one obtains a $B\sme A$-\ib{}
via the dual module $\X^\sim$: $\X^\sim$ coincides with $\X$ as a set
and, if $\iota:\X\to\X^\sim$ is the identity map, then the actions and
inner products are given by
\begin{alignat*}{2}
b\cdot\iota(x) &:= \iota(x\cdot b^*) &\qquad \blip
B<\iota(x),\iota(y)> &:= \rip B<x,y> \\
\iota(x)\cdot a &:= \iota(a^*\cdot x) & \qquad \brip
A<\iota(x),\iota(y)> &= \lip A<x,y>.
\end{alignat*}
Thus if $\bX$ is an $\A\sme\B$-\ibb, then there is a $\B\sme\A$-\ibb{}
$\bX^\sim$ which coincides with $\bX$ as a set and is such that the
identity map $\iota:\bX\to\bX^\sim$ maps the fibre $\X(u)$ onto
$\X(u)^\sim$.  Furthermore, $V_\g^\sim\(\iota(a)\):=\iota\(V_\g(a)\)$
defines a $G$-action on $\bX^\sim$ so that $(\bX^\sim,V^\sim)$
implements an equivalence between $(\B,\beta)$ and $(\A,\a)$.  Simply
put: ``$\sim$'' is symmetric.

Now suppose
\[
(\A,\a)\sim_{(\bX,V)}(\B,\beta)\sim_{(\bY,W)}(\CC,\gamma).
\]
Let $\X=\gsubogo(\bX)$ and $\Y=\gsubogo(\bY)$.
For each $u\in\go$, let $\ZZ(u):= \X(u)\tensor_{B(u)}\Y(u)$ be the
usual $A(u)\sme C(u)$-\ib{} (see, for example, \cite[p.~1036]{90c}).  
If $\set{\xi_1,\dots,\xi_n}\subseteq \X$ and
$\set{\eta_1,\dots,\eta_n} \subseteq \Y$, then
\[
\Bigl\| \sum_{i=1}^n \xi_i(u)\tensor \eta_i(u) \Bigr\|^2 = 
\Bigl\| 
\sum_{i,j=1}^n \blip \A<\xi_i(u), \xi_j(u) \cdot {\blip\B<\eta_j(u),
\eta_i(u)>} > \Bigr\|,
\]
and the right-hand side is continuous in $u$.  Moreover,
$\set{\xi(u)\tensor \eta(u):\text{$\xi\in\X$ and $\eta\in\Y$}}$ spans
a dense subset of $\ZZ(u)$ for each $u\in\go$.  Thus
\cite[Theorem~II.13.18]{fell-doran} implies that there is a Banach
bundle $\bZZ$ over $\go$ having $\ZZ(u)$ as fibres and such that
$u\mapsto \xi(u)\tensor\eta(u)$ belongs to $\ZZ=\gsubogo(\bZZ)$ for all
$\eta\in\Y$ and $\xi\in\X$.  Thus it remains only to see that the
$\A$- and $\B$-actions are continuous.  We consider only the
$\A$-action as the argument for the $\B$-action is similar.

Suppose that $(a_i,z_i)\to (a,z)$ in $\A*\bZZ$.  We aim to show that
$a_i\cdot z_i\to a\cdot z$.  For convenience, let $u=\pz(z)$ and
$u_i=\pz(z_i)$. 
Fix $\epsilon>0$ and choose $z_0:= \sum_{j=1}^n \xi_j(u)\tensor
\eta_j(u)$ such that each $\xi_j\in\X$, $\eta_j\in\Y$, and $\|a\cdot
z- a\cdot z_0\|<\epsilon$.  Therefore we eventually have
\[
\Bigl\| a_i\cdot \Bigl( \sum_{j=1}^n \xi_j(u_i) \tensor \eta_j(u_i)
\Bigr) - a_i\cdot z_i \Bigr\|<\epsilon,
\]
while 
\[
a_i \cdot \Bigl( \sum_{j=1}^n \xi_j(u_i) \tensor \eta_j(u_i)
\Bigr) \to a \cdot z_0.
\]
It follows from \cite[Proposition~II.13.12]{fell-doran} that $a_i\cdot
z_i\to a\cdot z$.  This completes the proof that ``$\sim$'' is
transitive, and hence that Morita equivalence of systems is an
equivalence relation.
\end{proof}

\begin{definition}
The collection $\brg$ of Morita equivalence classes of systems in $\sbrg$
is called the \emph{Brauer group of $G$}.
\end{definition}

Of course, the terminology anticipates that fact  
(\thmref{thm-grp}) that $\brg$ carries a
natural group structure that we shall define in a moment.

\begin{remark}
\label{rem-3.4}
Suppose that $(X,H)$ is a second countable locally compact transformation
group, and that $G=X\times H$ is the corresponding transformation group
groupoid.  Then our $\brg$ is the ``equivariant Brauer group''
$\br_H(X)$ of \cite{ckrw,krw,prw}.
\end{remark}

If $\A$ and $\B$ are elementary \cs-bundles over $T$, then there is
a unique elementary \cs-bundle $\A\tensor\B$ over $T\times T$ with
fibre $A(t)\tensor B(s)$ over $(t,s)$ and such that $(t,s)\mapsto
f(t)\tensor g(s)$ is a section for all $f\in A=\gsuboT(\A)$ and $g\in
B =\gsuboT(\B)$.
This bundle clearly satisfies Fell's condition if $\A$ and $\B$ do, as
does its restriction $\A\tensor_T\B$ to $\Delta=\set{(u,u)\in 
T\times T}$.
By identifying $\Delta$ with $T$, we will view $\A\tensor_T\B$ as an
elementary \cs-bundle over $T$.
Notice that 
$
\gsuboT(\A\tensor_T\B)
$
is the balanced tensor product $A\tensor_{C(T)}B$.
Consequently it should not prove confusing to write $A\tensor_T B$ for
$\gsuboT(\A\tensor_T\B)$. 

Now assume $(\A,\a),(\B,\beta)\in\sbrg$.  We want to equip $\batgob$
with a $G$-action $\alpha\tensor\beta$ such that
$(\batgob,\alpha\tensor\beta)\in\sbrg$.  Of course,
$\ag\tensor\beta_\g$ is a $*$-isomorphism of $A\(s_G(\g)\)\tensor B\(
s_G(\g) \)$ onto $A\(r_G(\g)\)\tensor B\(
r_G(\g) \)$; thus, $\batgob$ admits a $G$-action provided we show that
this action is continuous.
There is a map $(a,b)\mapsto a\tensor b$ from
$\A*\B=\set{(a,b)\in\A\times \B: \pa(a)=\pb(b)}$ into $\batgob$, and
it is clearly continuous.  
Suppose that $(\gi,c_i)\to (\g,c)$ in $\(\batgob\)*G$.  Let
$u=s_G(\g)$ and $u_i=s_G(\gi)$.  Fix $\epsilon>0$.
Choose $z=\sum_{j=1}^n a^j\tensor b^j\in \atgob(u)$ such that
$\|z-c\|<\epsilon$.  Choosing sections in $A$ and $B$ through $a^j$
and $b^j$, respectively, we obtain elements $(a_i^j,b_i^j)\to
(a^j,b^j)$ in $\A*\B$ with $\pa(a^j_i)=u_i=\pb(b^j_i)$.  Then 
\[
z_i:= \sum_{j=1}^n a_i^j\tensor b_i^j\to z.
\]
On the other hand, $\gi\cdot z_i\to \g\cdot z$.  Moreover, $\|\g\cdot
c - \gamma\cdot z\|=\|c-z\|<\epsilon$, and $\|\gi\cdot z_i - \gi\cdot
c_i\| = \|z_i - c_i\|$.  Therefore we eventually have $\|\gi\cdot z_i
- \gi\cdot c_i\| < \epsilon$.  It follows that $\gi\cdot c_i\to
\g\cdot c$ as required.

Now suppose that 
\begin{gather*}
\aa\sim_{\xv}(\CC,\g)\text{, and}\\
\bb\sim_{(\bY,W)}(\D,\delta).
\end{gather*}
For each $u\in\go$, the external tensor product
$\ZZ(u):=\X(u)\tensor\Y(u)$ is an $A(u)\tensor B(u)\sme C(u)\tensor
D(u)$-\ib{} (see \cite[\S1.2]{je-th} or \cite[\S2]{ckrw}). 
If $\xi_1,\dots,\xi_n\in\X$ and $\eta_1,\dots,\eta_n\in\Y$, then
\[
\Bigl\|\sum_{i=1}^n \xi_i(u) \tensor\eta_i(u) \Bigr\|^2 = \Bigl\|
\sum_{i,j=1}^n \brip\CC <\xi_i(u),\xi_j(u)>\tensor \brip \D <\eta_i(u)
, \eta_j(u)> \Bigr\|,
\]
and the right-hand side is continuous in $u$.  Since $\set{f(u)\tensor
g(u):\text{$f\in\X$ and $g\in\Y$}}$ spans a dense subspace of $\ZZ(u)$,
there is a Banach bundle $\bZZ:=\bX\tensor\bY$, having $\ZZ(u)$ as fibre
over $u$, such that $u\mapsto f(u)\tensor g(u)$ is a section for
all $f\in\X$ and $g\in\Y$.  Note that $\bZZ$ is the bundle associated
to the global external tensor product $\X\tensor\Y$.  There are
isometries $V_\gamma\tensor W_\g:\X\(s_G(\g)\)\tensor \Y\(s_G(\g)\)
\to \X\(r_G(\g)\)\tensor \Y\(r_G(\g)\)$ which define a left $G$-action
by ismorphisms on $\bZZ$, and it is not hard to see that this action is
continuous.  Thus $(\bX\tensor\bY,V\tensor W)$ implements an
equivalence between $(\batgob,\a\tensor\beta)$ and
$(\CC\tensor_{\go}\D,\g\tensor \delta)$.

\begin{definition}
Let $\TT_{\go}=\TT=\go\times\C$ be the trivial line bundle over $\go$,
and let 
$\tau$ be the action of $G$ on $\TT$ given by $\g\cdot\(s_G(\g),\lambda):=
\(r_G(\g),\lambda)$.
\end{definition}

\begin{prop}
The binary operation
\begin{equation}
\label{eq-mult}
[\A,\a][\B,\beta]:=[\batgob,\alpha\tensor\beta]
\end{equation}
is well-defined on $\brg$.  With respect to this operation, $\brg$ is
an abelian semigroup with identity equal to the class
$[\TT,\tau]$. 
\end{prop}

\begin{proof}
The preceding discussion shows that the operation is well-defined.
Since equivariant Banach bundle isomorphisms certainly provide Morita
equivalences of systems, it is not hard to check that the operation is
commutative and associative.  

The isomorphism $\phi\tensor f\mapsto \phi\cdot f$ of
$C_0(\go)\tensor_{\go} A$ onto $A$ induces a Banach bundle isomorphism
of $\TT\tensor_{\go}\A$ onto $\A$ taking $\lambda\tensor a$ to
$\lambda a$.  Since this map is equivariant, it follows that
$[\A,\a][\TT,\tau]= [\A,\a]$; that is, $[\TT,\tau]$ is an identity as
claimed. 
\end{proof}

In order to see that $\brg$ is a group, we will need to introduce the
\emph{conjugate bundle}.  If $p:\A\to\go$ is a Banach bundle, then let
$\bA$ be the topological space $\A$ and let $\iota:\A\to\bA$ be the
identity map.  Then $\bp:\bA\to\go$ defined by $\bp\(\iota(a)\)=
\iota\(p(a)\)$ is a Banach bundle over $\go$ with fibre
$\overline{A}(t)$
identified
with the conjugate Banach space $\overline{A(t)}$ (see the discussion
preceding \cite[Remark~3.5]{ckrw}).  Furthermore
$\gsubogo(\bA)=\overline{\gsubogo(\A)}$, so we may write
$\overline{A}$ for the section algebra of $\bA$.  If $\aa\in\sbrg$,
then, $\ba_\g\(\iota(a)\):= \iota\(\ag(a)\)$ is a \cs-$G$-action and
$(\bA,\ba)\in \sbrg$.

\begin{thm}
\label{thm-grp}
If $G$ is a locally compact groupoid with paracompact unit space
$\go$, then $\brg$ is an abelian
 group with addition defined by \eqref{eq-mult}.  
The inverse of $[\A,\a]$ is $[\bA,\ba]$.
\end{thm}

\begin{proof}
At this point, it will suffice to show that for any $\aa\in\sbrg$,
$(\A\tensor_{\go}\bA,\a\tensor\ba)\sim(\TT,\tau)$.  Our proof will
follow the lines of \cite[Theorem~3.6]{ckrw}, and we will use some of
the same notation.  In particular, since $A=\gsubogo(\A)$ has
continuous trace, it follows from \cite[Lemmas 6.1~and 6.2]{90c} that
there are compact sets $F_i\subseteq \go$ whose interiors form a cover
$\U=\set{\operatorname{int}(F_i):i\in I}$ of $\go$ such that
\begin{enumerate}
\item
for each $i\in I$, there is an $A^\Fi\smeover\Fi C(\Fi)$-\ib{} $\X_i$,
and
\item
for each $i,j\in I$, there is an \ib{} isomorphism\\
$\gij:\X_j^\fij\to\X_i^\fij$. 
\end{enumerate}
Consequently, there is an \ibb{} $\bX_i$ with $\gsubogo(\bX_i)=\X_i$, and
the fibres $\X_i(u)$ are Hilbert spaces.  If $u\in\go$, then $\gij$
determines a Hilbert space isomorphism $\gij(t):\X_j(t)\to \X_i(t)$.
The \emph{Dixmier-Douady class} of $\A$ is determined by the cocycle
$\nu=\set{\nu_{ijk}}\in H^2(\U,\mathcal S)$, where for each $t\in
F_{ijk}$
\[
\gij(t)\circ g_{jk}(t)=\nu_{ijk}(t) g_{ik}(t).
\]
We can then define $\bar\gij:\oX_j^\fij\to\oX_i^\fij$ by
$\bar\gij\(\iota(x)\) := \iota\(\gij(x)\)$, and $h_{ij}:= \gij\tensor
\bar\gij : \X_j^\fij \tensor_{C(\fij)} \oX_j^\fij \to \X_i^\fij
\tensor_{C(\fij)} \oX_i^\fij$ as in \cite{ckrw}.  We obtain an
$A\tensor_{\go} \overline{A}\smeover\go C_0(\go)$-\ib{} by forming
\[
\Y':=\set{(y_i)\in\prod\X_i\tensor_{C(\Fi)}\oX_i: h_{ij}(t)\(y_j(t)\)=y_i(t)},
\]
and noticing that if $t\in\fij$ and $(x_k),(y_k)\in\Y'$, then
\[
\rip C(\Fi)<x_i,y_i>(t)=\rip C(F_j)<x_j,y_j>(t).
\]
Then we can define $\brip C_0(\go)<(x_k),(y_k)>(t) = \rip
C(\Fi)<x_i,y_i>(t)$ when $t\in\Fi$, and
\[
\Y=\set{y\in\Y':\text{$t\mapsto \rip C_0(\go)<y,y>(t)$ vanishes at
$\infty$}}. 
\]
It is shown in \cite{ckrw,90c} that $\Y$ is an
$A\tensor_{\go}\overline{A} \sme C_0(\go)$-\ib, and that $\Y$ is
isomorphic as a (full) Hilbert $C_0(\go)$-module to
\[
\NN =\set{a\in A:\text{$t\mapsto \tr\(a^*a(t)\)$ is in $C_0(\go)$}}.
\]
The isomorphism $\Phi:\Y\to\NN$ is given by $\Phi\((y_k)\)(t) =
\Phi_i(y_i)(t)$ if $t\in\Fi$, where
$\Phi_i:\X_i\tensor_{C(\Fi)}\oX_i\to\NN^\Fi$ is given by
$\Phi_i\(x\tensor\iota(y)\)(t) =\lip A^\Fi<x,y>(t)$.  Now
$\NN=\gsubogo(\bNN)$ for some Banach bundle $\bNN$.  It also follows
from the uniqueness condition in \cite[Theorem~II.13.18]{fell-doran}
that we can identify $\bNN$ with
\[
\set{a\in \A:\tr(a^*a)<\infty}
\]
\emph{with the relative topology}.  (The fibres of $\bNN$ are the
ideals of
Hilbert-Schmidt operators in the appropriate fibre of $\A$.)
Since each $\ag$ is an isomorphism
of elementary \cs-algebras, and therefore trace preserving,
$\ag(\bNN)\subseteq \bNN$.  It follows that we can define an
continuous action $V$ of
$G$ on $\bNN$ by isomorphisms via $V_\g(a):=\ag(a)$.  

\begin{lem}
\label{lem-ckrw3.14}
Let $(\B,\tau)\in\sbrg$.
Suppose that $q:\bX\to T$ is an $\A\sme\B$-\ibb, and that $V$ is a
continuous action of $G$ on $\bX$ by isomorphisms such that
$\tau_\g\(\rip\B<x,y>\cdot b\) =
\brip\B<V_\g(x),V_\g(y)>\cdot\tau_\g(b)$.
Then $\alpha_\g(a):= V_\g a V_\g^{-1}$ defines an isomorphism of
$A\(s_G(\g)\) $ onto $A\(r_G(\g)\)$, such that $\aa\in\sbrg$, and
$\aa\sim_{\xv}(\B,\tau)$. 
\end{lem}
\begin{proof}
Each $a\in A(u)$ defines a linear operator $x\mapsto a\cdot x$ in
$\L\(\X(u)\)$, and we can identify $A(u)$ with $\K\(\X(u)\)$.  If
$\g\in G$ and $a\in A\(s_G(\g)\)$, then $x\mapsto V_\g aV_\g^{-1}$ is
clearly in $\L\(\X\(r_G(\g)\)\)$.  If $z\in \X\(r_G(\g)\)$, then
\begin{align*}
V_\g\lip\A<x,y>V_\g^{-1}(z) 
&= V_\g\( x \cdot \brip\B<y,V_\g^{-1}(z)>\) \\
&= V_\g(x) \cdot \tau_\g\(\brip\B<y,V_\g^{-1}(z)>\) \\
&= V_\g(x) \cdot \brip\B<V_\g(x),z> \\
&= \brip\A<V_\g(x),V_g(y)>\cdot z.
\end{align*}
It follows that $V_\g\lip\A<x,y>V_\g^{-1}\in \K\(\X\(r_G(\g)\)\)\cong
A\(r_G(\g) \)$.
Thus we can view $V_\g a V_\g^{-1}$ as an element of $A\(r_G(\g)\)$ for
each $a\in A\(s_G(\g)\)$.  The continuity of $V$ ensures that $\alpha$
is continuous.  The rest is routine.
\end{proof}\renewcommand{\qed}{}\end{proof}

\begin{proof}
[Completion of the proof of \thmref{thm-grp}]
Note that
$\K(\NN)$ is Morita equivalent to $C_0(\go)$, and, since $\Y$ and
$\NN$ are isomorphic, there is an isomorphism
$Q:A\tensor_{\go}\overline{A}\to\K(\NN)$ defined by
$Q(a\tensor\iota(b))(n) = acb^*$ (\cite[Corollary~3.12]{ckrw}).
It is clear from the uniqueness of the trace on an elementary
\cs-algebra that 
\[
\brip C_0(\go)<V_\g(n),V_\g(m)>\(r_G(\g)\) = \rip
C_0(\go)<n,m>\(s_G(\g)\). 
\]
Therefore it follows from \lemref{lem-ckrw3.14} and the fact that
$\K(\NN)$ is Morita equivalent to $C_0(\go)$, that there is a
\cs-$G$-bundle action $\beta$ on $\A\tensor_{\go}\bA$ such that
\[
\(\A\tensor_{\go}\bA,\beta)\sim(\TT,\tau).
\]
In particular
\begin{align*}
Q\circ \beta_\g\(f(\g)\tensor \iota\(g(\g)\)\) \cdot c 
&= V_\g Q\(f(\g)\tensor\iota\(g(\g)\)\) V_\g^{-1}(c) \\
&= \alpha_\g\(f(\g)\alpha_\g^{-1}(c)g(\g)^*\) \\
&= Q\(\ag\(f(\g)\) \tensor \iota\(\ag\(g(\g)\)\)\)(c).
\end{align*}
Therefore $\beta=\alpha\tensor\bar\alpha$.  This completes the proof
that $\brg$ is a group.
\end{proof}
\else

\fi
\ifdraft
%
%

\section{The First Isomorphism Result}\label{sec-fiso}

The main objective of this section is to prove that equivalent
groupoids have isomorphic Brauer groups.  Specifically:

\begin{thm}
\label{thm-main1}
Suppose that $Z$ is a \equi(H,G).  Then
$\phi^Z\([\A,\a]\):=[\A^Z,\alpha^Z]$ defines a group
isomorphism of
$\brg$ onto $\br(H)$.
\end{thm}

This is the
natural groupoid analogue of the main result in \cite{krw} (see
\propref{equivalent}).

The main tool here will be the fundamental construction of
\propref{prop-fundcon} (which applies to any Banach bundle).  
Our first result is that this construction
respects Morita equivalence so that $\phi^X$ is well defined.

\begin{lem}
\label{lem-wd}
Suppose that $Z$ is a \equi(H,G), and that $\aa\sim_\xv(\B,\beta)$ in
$\sbrg$.
Then $(\A^Z,\alpha^Z)\sim_{(\bX^Z,V^Z)} (\B^Z,\beta^Z)$ in $\sbr(H)$.
\end{lem}

\begin{proof}
As usual, we will denote the section algebras of $\A^Z$, $\B^Z$, and $\bX^Z$
by $A^Z$, $B^Z$, and $\X^Z$, respectively.    
We obtain a well defined continuous map from $\A^Z*\bX^Z$ to
$\bX^Z$ by
\[
[z\cdot\g,a]\cdot[z\cdot\eta,\xi]:= [z,\ag(a)\cdot V_\eta(\xi)];
\]
similarly, there is a continuous map from $\bX^Z*\B^Z$ to $\bX^Z$
given by
\[
[z\cdot\eta,\xi] \cdot [z\cdot\gamma,b] = [z,V_\eta(\xi)\cdot\beta_\g(b)].
\]
Furthermore, there are well defined maps $\rip\B^Z<\cdot,\cdot>$ and
$\lip\A^Z<\cdot, \cdot>$ from $\bX^Z*\bX^Z$ to $\B^Z$ and $\A^Z$,
respectively, given by
\begin{align*}
\brip\B^Z<[z\cdot\g,\xi],[z\cdot\eta,\xi']> &=
\bigl[z,\rip\B<\xi,\xi'>\bigr] \text{, and} \\
\blip\A^Z<[z\cdot\g,\xi],[z\cdot\eta,\xi']> &=
\bigl[z,\lip\A<\xi,\xi'>\bigr] .
\end{align*}
It is not hard to see that these formulae equip
$\X^Z(v)$ with an $A^Z(v)\sme B^Z(v)$-\ib{} structure.  Since the norm
on $\X^Z(v)$ induced by the inner products coincides with the Banach
bundle norm, it follows that $\bX^Z$ is an $\A^Z\sme \B^Z$-\ibb.
Since
\begin{align*}
\brip\B^Z<V_h^Z\([z\cdot\g,\xi]\), V_h^Z\([z\cdot\eta,\xi']\)> &=
\brip\B^Z <[h\cdot z\cdot\g,\xi],[h\cdot z\cdot\eta,\xi']> \\
&= \bigl[h\cdot z , \brip\B<V_\g(\xi),V_\eta(\xi')> \bigr] \\
&= \beta_h^Z\(\bigl[z,  \brip\B<V_\g(\xi),V_\eta(\xi')> \bigr] \) \\
&= \beta_h^Z\( \brip\B^Z < [z\cdot\g,\xi],[z\cdot\eta,\xi']> \),
\end{align*}
and similarly for $\alpha^Z$, we see that $(\bX^Z,V^Z)$ implements
the desired equivalence.
\end{proof}

The next lemma will be exactly what we need to prove that $\phi^X$ is
a homomorphism.

\begin{lem}
Suppose that $\pa:\A\to T$ and $\pb:\B\to T$ are elementary
\cs-bundles, and that $q:X\to T$ is continuous.  Then the map
$(x,\sum_{i=1}^n a_i\tensor b_i)\mapsto \sum_{i=1}^n (x,a_i)\tensor
(x,b_i)$ extends to a Banach bundle isomorphism of $q^*\(\A\tensor_T\B)$ onto
$q^*\A\tensor_X q^*\B$.
\end{lem}

\begin{proof}
The given map is isometric on fibres and therefore clearly extends to
a map $\Phi$ from $q^*\(\A\tensor_T\B)$ onto
$q^*\A\tensor_X q^*\B$ which preserves fibres and is an isometry on
each fibre.  Thus it will follow that $\Phi$ is a Banach bundle
isomorphism once we see that it is continuous.  But if $f$ and $g$ are
sections of $\A$ and $\B$, respectively, then $q^*(f\tensor g)$ given
by $x\mapsto \(x,f\(q(x)\)\tensor g\(q(x)\)\)$ is a section of
$q^*(\A\tensor_T\B)$ and $\Phi\circ q^*(f\tensor g) = q^*f\tensor
q^*g$ is a section of $q^*\A\tensor_X q^*\B$ (here $q^*f(x) :=
\(x,f\(q(x)\)\)$ and $q^*g(x):= \(x,g\(q(x)\)\)$).  Now the result
follows from \cite[Proposition~II.13.16]{fell-doran}.
\end{proof}

\begin{cor}
\label{cor-hm}
Suppose that $\aa,(\B,\beta)\in\sbrg$ and that $X$ is a \equi(H,G).
Then $\((\atgob)^X,(\alpha\tensor\beta)^X\)$ is equivalent to
$\(\A^X\tensor_{\ho} \B^X, \a^X\tensor\beta^X\)$ in $\sbr(H)$.
\end{cor}

\begin{proof}
The result follows immediately from the lemma and the observation that
$\bigl[(x,a)\tensor (x,b)\bigr]\mapsto [x,a]\tensor [x,b]$ is an
$H$-invariant Banach bundle isomorphism of
$\(\sx^*\A\tensor_X\sx^*\B\)/G$ onto $\A^X\tensor_{X/G}\B^X$.
\end{proof}

It follows from \lemref{lem-wd} and \corref{cor-hm} that if $X$ is a
\equi(H,G), then $\phi^X$ is a
well defined homomorphism from $\brg$ into $\br(H)$.  We can finish
the proof of \thmref{thm-main1} by exhibiting an inverse for
$\phi^X$.  Our next result will provide the necessary ``calculus'' of
these maps; this will allow us not only to finish  the proof of
\thmref{thm-main1}, but will be useful in \secref{sec-e(x)} as well.

First recall that groupoid equivalence is indeed an equivalence
relation.  Reflexivity follows from noting that $G$ is a \equi(G,G).  If
$X$ is a \equi(H,G), then the opposite space $\xop$ is a \equi(G,H): if
$\iota:X\to\xop$ is the identity map then the left $G$-action on
$\xop$ is given by $\g\cdot\iota(x):=\iota(x\cdot\g^{-1})$ and the
right $H$-action is given by $\iota(x)\cdot h:=\iota(h^{-1}\cdot x)$.
If $X$ is a \equi(H,G) and $Y$ is a \equi(G,K), then we obtain an \equi(H,K)
by forming $\xsgy$ which is the orbit space of $X*Y=\set{(x,y)\in
X\times Y: \sx(x)=r_Y(y)}$ by the diagonal $G$ action:
$(x,y)\cdot\g:=(x\cdot \g,\g^{-1}\cdot y)$. 
Furthermore two $(G_1,G_2)$-equivalences $Z_1$ and $Z_2$ are
\emph{isomorphic} if there is a homeomorphism $\phi:Z_1\to Z_2$ such
that $\phi(\g_1\cdot z\cdot\g_2)=
\g_1\cdot \phi(z)\cdot \g_2$ for all $(\g_1,z,\g_2)\in G_1*Z*G_2$.

 More generally, if $X$ is
a principal right $G$-space, then we will also write $\xop$ for the
corresponding principal left $G$-space.  Then we can form $G^X:=\imp(X)$ as
above.   Then $G^X$ has a natural groupoid structure (cf., e.g.,
\cite[pp.~119+]{gcddi}).  In particular, the unit space of $G^X$,
$\set{[x,x]:x\in X}$, can be identified with $X/G$, and $G^X$ acts on
the left of $X$ via $[y,x]\cdot x=y$.  With this action $X$ is a
\equi(G^X,G), and we refer to $G^X$ as the \emph{imprimitivity
groupoid} of $X$ (\cite[Theorem~3.5(1)]{gcddi}). 

\begin{lem}
\label{lem-inv}
Suppose that $G$, $H$, and $K$ are locally compact groupoids.
\begin{enumerate}
\item
If $X$ is a \equi(H,G) and $Y$ is a \equi(G,K), then
\[
\phi^{\xsgy}=\phi^X\circ\phi^Y.
\]
\item
If $X$ and $Y$ are isomorphic as \equi(H,G)s, then $\phi^X=\phi^Y$.
\item
Viewing $G$ is a \equi(G,G), $\phi^G=\id_G$.
\item
If $X$ is a \equi(H,G), then
\[
\phi^{\xopshx}=\id_{\brg}\quad\text{and}\quad\phi^{\xsgxop}=\id_{\br(H)}.
\]
\item
If $X$ is a \equi(H,G) and $Y$ is a \equi(K,G), then
\[
\phi^X= \phi^Y\circ\phi^{X\understar{G}\opposite{Y}}.
\]
\end{enumerate}
\end{lem}

\begin{proof}
Suppose that $\aa\in\sbr(K)$.  We need to prove that
$(\A^\xsgy,\a^\xsgy) \sim \((A^Y)^X,(\a^Y)^X\)$.  But
$\bigl[x,[y,a]\bigr] = \bigl[x_1,[y_1,a_1]\bigr]$ in 
\[
(\A^Y)^X=\set{\bigl[x,[y,a]\bigr]:\text{$\sx(x)=r_Y(y)$ and
$s_Y(y)=\pa(a)$}}
\]
if and only if there are $\g,\eta\in G$ such that $x_1=x\cdot\g$,
$y_1= \g\cdot y\cdot \eta$, and $a_1=\a_\eta^{-1}(a)$.  That is, if
and only if $\bigl[[x,y],a\bigr] = \bigl[[x_1,y_1],a_1\bigr]$ in
\[
\A^\xsgy =\set{\bigl[[x,y],a\bigr]: \text{$\sx(x)=r_Y(y)$ and
$s_Y(y)=\pa(a)$}}.
\]
Thus we have a well defined bijection $\bigl[x,[y,a]\bigr] \mapsto
\bigl[[x,y],a\bigr] $ from $(\A^Y)^X$ onto $\A^\xsgy$.  It is not hard
to see that this is an equivariant Banach bundle map.  Part~\partref1
follows.

Parts \partref2 and \partref3 are straightforward.  Part~\partref4
follows from \partref2~and \partref3.  (It should be noted that $H$ is
only \emph{isomorphic} to $\xsgxop$ and one may wonder if this
isomorphism intervenes in \partref4.  The point of \partref2 is that
it doesn't.)  Finally, \partref5 follows from
\partref3.
\end{proof}

\begin{proof}
[Proof of \thmref{thm-main1}]
We have already observed that $\phi^X$ is a \hm.
It follows from \lemref{lem-inv} that $\phi^{\xop}$ is an
inverse for $\phi^X$.
\end{proof}
\else

\fi
\ifdraft
%
%

\section{Haar Systems for Imprimitivity Groupoids}\label{sec-haar}

Suppose that $X$ and $Y$ are locally compact (Hausdorff) spaces and
that $\pi:X\to Y$ is a continuous, open, surjection.  Then $C_c(X)$ is
a $C_c(Y)$-module.  A $C_c(Y)$-module map $\beta:C_c(X)\to C_c(Y)$
with the property that $f\ge0$ implies $\beta(f)\ge0$ is called a
\emph{$\pi$-system}.  Notice that for each $y\in Y$ there is a Radon
measure
$\beta_y$ with $\supp(\beta_y)\subseteq \pi^{-1}(y)$ such that
\[
\beta(f)(y)=\int_X f(x)\,d\beta_y(x).
\]
We say that $\beta$ is \emph{full} if $\supp(\beta_y)=\pi^{-1}(y)$ for
all $y\in Y$.  If in addition, $X$ and $Y$ are (right) $G$-spaces,
then we say that $\beta$ is \emph{equivariant} if 
\[
\int_X f(x\cdot \g)\,d\beta_y(x) = \int_X f(x)\,d\beta_{y\cdot
\g}(x)\quad
\text{for all $(y,\g)\in Y*G$.}
\]

We will make use of the following result from \cite{renault2} where it
was proved for left actions.  We restate it for convenience.

\begin{lem}
[{\cite[Lemme 1.3]{renault2}}]%
\label{lem-1.3}
Suppose that $X$ and $Y$ are principal right $G$-spaces, and that
$\pi:X\to Y$ is a continuous, open, equivariant surjection.  Let
$\dot\pi:X/G \to Y/G$ be the induced map.
\begin{enumerate}
\item
If $\beta$ is an equivariant $\pi$-system, then
\[
\dot\beta(f)(y\cdot G):= \int_{X/G} f(x\cdot G)\,d\beta_y(x)
\]
is a well-defined $\dot\pi$-system.
\item
If $\tau$ is a $\dot\pi$-system, then there is a unique equivariant
$\pi$-system $\beta$ such that $\dot\beta=\tau$.
\end{enumerate}
\end{lem}

Our interest in \lemref{lem-1.3} here is that it provides a
characterization of when imprimitivity groupoids carry a Haar system.

\begin{prop}
\label{prop-haar}
Suppose that $G$ and $H$ are locally compact groupoids and that $X$ is
a \equi(H,G).  Then $H$ has a Haar system if and only if $X$ has a
full equivariant $s$-system where $s:X\to\go$ is the source map for
the principal $G$ action on $X$.
\end{prop}

\begin{proof}
We can identify $H$ with $G^X=\xsgxop$.  Suppose that $\alpha$ is a
full equivariant $s$-system on $X$.  Let $\pil:X*\xop\to X$ be the
projection map: $\pil\(x,\iota(y)\)=x$.  Let $\beta_x:=\delta_{\set
x}\times \alpha_{s(x)}$.  Then it is not hard to verify that
\[
\beta(f)(x)=\int_{X*\xop} f\(x,\iota(y)\)\,d\alpha_{s(x)}(y)
\]
is a full $\pil$-system.
Furthermore since $\alpha$ is equivariant,
\begin{align*}
\int_{X*\xop} f\(x,\iota(y)\)\,d\beta_{x\cdot\g}(x,y) & =
\int_{X*\xop} f\(x\cdot \g,\iota(y)\)\,d\alpha_{s(x\cdot\g)}(y) \\
&= \int_{X*\xop} f\(x\cdot\g,\iota(y\cdot\g)\,d\alpha_{s(x)}(y) \\
&= \int_{X*\xop} f\(x\cdot\g,\g^{-1}\cdot\iota(y)\)\,d\beta_x(x,y);
\end{align*}
it follows that $\beta$ is equivariant as well.  Using the first part
of \lemref{lem-1.3}, we obtain a $\dot\pil$-system
$\dot\beta: C_c(G^X)\to C_c(X/G)$.  It is not hard to verify that
$\l^{x\cdot G}:=\beta_{x\cdot G}$ is a Haar system for $G^X$.  For
example, invariance follows from the  computation:
\begin{align*}
\int_{G^X} f\([x,y][u,v]\)\,d\l^{y\cdot G}\([u,v]\)
&= \int_{G^X} f\([x,v]\) \,d\beta_{y\cdot G}
=  \int_{G^X} f\([x,v]\) \,d\alpha_{s(y)}(v) \\
&= \int_{G^X} f\([x,v]\) \,d\alpha_{s(x)}(v)\text{, since
$s(x)=s(y)$,} \\
&= \int_{G^X} f\([x,v]\) \,d\l^{x\cdot G}\([x,v]\).
\end{align*}

For the converse, notice that a Haar system $\l$ for $G^X$ defines a
$\dot\pil$-system.  Using the second part of \lemref{lem-1.3}, there is a
unique $\pil$-system $\beta$ such that $\dot\beta$ coincides with
$\l$.  But $\beta$ must be of the form $b_x=\delta_{\set
x}\times\alpha_x$ with $\supp{\alpha_x}=s^{-1}(x)$.  Furthermore,
$\alpha_x$ depends only on $s(x)$ and
$\alpha=\set{\alpha_u}_{u\in\go}$ is a full equivariant $s$-system as
required. 
\end{proof}

In the sequel, we write $\pg$ for the collection of second
countable\footnote{See the footnote on page~\pageref{sec-four}.}
principal right $G$-spaces which have full equivariant $s$-systems.
By \propref{prop-haar}, these are precisely the principal $G$-spaces
whose imprimitivity groupoids have Haar systems.  We shall see in
Sections \ref{sec-graph}~and \ref{sec-four} that $\pg$ may be viewed
as a generalization of the collection of all open covers of $\go$,
directed by refinement.

\else

\fi
\ifdraft
%
%
\section{The Graph of a Homomorphism}\label{sec-graph}

This section contains two key results for our analysis.  The first,
\lemref{lem-tech}, gives a serviceable sufficient condition for a \hm{}
between two groupoids to implement an equivalence.  The second,
\corref{cor-trivialdd}, shows that if $[\A,\a]\in\brg$, then there is
a space $X\in\pg$ such that $\delta(\A^X)=0$.  Our approach to
understanding $\brg$ is to analyze it ``through the eyes'' of $\A^X$.

If $\phi:H\to G$ is any continuous groupoid \hm{}, then we can define
a principal right $G$-space called the \emph{graph of $\phi$}
as follows:
\[
\Gr(\phi):= \set{(u,\gamma)\in \ho\times G:\phi(u)=r(\gamma)}.
\]
If, for example, $\phi\restr\ho$ is an open surjection onto $\go$, then 
$s(u,\g):=s(\g)$ is clearly an open surjection of $\grp$ onto
$\go$, and $\grp$ is a principal right $G$-space:
$(u,\gamma)\cdot\g'= (u,\g\g')$.
Moreover, we claim that $r(u,\g):=r(\g)$ is open\footnote{We have
opted not to decorate every range and source map in the sequel.
In particular, writing
$r_{\grp}(u, \g)=r_G(\g)$ in place of $r(u,\g):=r(\g)$ does not seem
to make these expressions easier to sort out.  Instead, 
we hope our meaning will be clear from context.}.
If $v_i\to v$ in $\go$, then we can \psnrl, and assume that there are
$\gi\to\g$ in $G$ 
such that $r(\g_i)=v_i$.  Since $\phi\restr\ho$ is open and
surjective, we can \psnrl{} and assume that there are $u_i\to u$ 
in $\ho$ and $\phi(u_i)=v_i$.  Thus $(u_i,\g_i)$
converges to $(u,\g)$ in $\grp$; it follows that $s$ is open.
Therefore $\grp$ is also a left $H$-space: 
$h\cdot\(s(h),\g\)=\(r(h),\phi(h)\g\)$.   We are most
interested in situations where $\grp$ is actually a
$(H,G)$-equivalence; in this case we say that $\phi$
\emph{induces an equivalence between
$H$ and
$G$}. 

\begin{ex}\label{ex-4.1}
Suppose that $\phi$ is a groupoid automorphism of $G$.  Then $\grp$ is
easily seen to be a $(G,G)$-equivalence, and $\phi$ induces an
equivalence between $G$ and $G$.
\end{ex}

Our next lemma gives some simple criteria for a \hm{} to induce an
equivalence; we have made no attempt to find the most general such
criteria --- only ones that are easy to check in our applications.

\begin{lem}
\label{lem-tech}
Suppose that $\phi:H\to G$ is a continuous groupoid \hm{} that
satisfies the following.
\begin{enumerate}
\item
$\phi$ is surjective.
\item
$\phi\restr{\ho}$ is open.
\item
If $\set{h_i}$ is a net in $H$ such that $\set{r(h_i)}$ and
$\set{s(h_i)}$ converge in $\ho$, while $\set{\phi(h_i)}$ converges in
$G$, then $\set{h_i}$ has a convergent subnet.
\item
The restriction of $\phi$ to the isotropy subgroupoid $A=\set{h\in
H:r(h)=s(h)}$ has  trivial kernel; that is, $\ker(\phi)\cap A=\ho$, where
$\ker(\phi)=\phi^{-1}(\go)$. 
\item
If $u,v\in\ho$ and if $\phi(u)=\phi(v)$, then there exists
$h\in\ker(\phi)$ such that $h\cdot v=u$.
\end{enumerate}
Then $\phi$ induces an equivalence between $H$ and $G$.  If only
assertions \partref1~and
\partref2 hold, then $\grp$ is a principal right $G$-space and a
$H$-space such that $h\cdot(x\cdot\g)=(h\cdot x)\cdot\g$ and
$\grp/G\cong \ho$.
\begin{proof}
The $G$-action on $\grp$ is always free and proper (provided that
$\phi(\ho)=\go$ and $\phi\restr\ho$ is open).  If $\phi\restr\ho$ is
open, then, as pointed out above, $\grp$ is a right $H$-space, and
this action will be free and proper in view of \partref4 and
\partref3, respectively.  It remains to show only that $s:\grp\to\go$
and $r:\grp\to\ho$ induce homeomorphisms of $\go$ with $\hmgrp$ and
$\ho$ with $\grpmg$, respectively.

Assume that $(u,\g)$ and $(v,\eta)$ are in $\grp$ and that
$s(u,\g)=s(\g)=s(\eta)=s(v, \eta)$.  By \partref1, there is an $h\in
H$ such that $\phi(h)=\eta\g^{-1}$.  Then $\phi\(s(h)\) =
r(\g)=\phi(u)$ and $\phi\(r(h)\)=r(\eta)=\phi(v)$.  It follows from
\partref5, that there are $h',h''\in\ker(\phi)$ such that $h'\cdot
r(h)=v$ and $h''\cdot u=s(h)$.  Then $\phi(h'hh'')=\eta\g^{-1}$ and
$h'hh''\cdot (u,\g)=(v,\eta)$.  Thus $s$ induces a bijection of
$\hmgrp$ onto $\go$.  Now suppose that $w_i\to w$ in $\go$ and that
$(u,\g)\in \grp$ satisfies $s(u,\g)=s(\g)=w$.  Passing to a subnet and
relabeling we can assume that there are $\g_i\in G$ such that
$\g_i\to\g$ and $s(\g_i)=w_i$.  In view of \partref2, we can
even assume that there are $u_i\to u$ with $\phi(u_i)=r(\g_i)$.  Then
$(u_i,\g_i)\to (u,\g)$ in $\grp$.  It follows that $s$ is open and
induces a homeomorphism of $\hmgrp$ with $\go$ as claimed.

The case for $r:\grp\to\ho$ is more straightforward.  It is immediate
that $r$ induces a continuous bijection of $\grpmg$ onto $\ho$.
Furthermore if $(u,\g)\in\grp$ and if $u_i\to u$, then we can \psnrl{}
and assume there are $\g_i\to \g$ in $G$ with $r(\g_i)=\phi(u_i)$.  Then
$(u_i,\g_i)\in\grp$, and $(u_i,\g_i)\to (u,\g)$.  Therefore $r$ is
open and always induces a homeomorphism of $\grpmg$ onto $\ho$.
\end{proof}
\end{lem}

\begin{ex}
[Local Trivialization]
\label{ex-loctriv}
Let $\U=\set{U_i}$ be a locally finite open cover of $\go$. 
Let
$X:=\coprod U_i$ and $\psi:X\to\go$ the obvious map --- note that
$\psi$ is a continuous, open surjection.  Then we can build a groupoid
$G^\psi:=\coprod_{ij} G_{U_i}^{U_j}$ as follows.  We write
elements of
$G^\psi$ as triples $(j,\g,i)$ with $r(\g)\in U_j$ and
$s(\g)\in U_i$.  Each
$(i,u)\in X$ is identified with $(i,u,i)\in G^\psi$.  Thus we can
define $r(j,\g,i)=\(j,r(\g)\)$ and $s(j,\g,i)=\(i,s(\g)\)$.  Then we
can identify the unit space of $G^\psi$ with $X$ where the
groupoid operations on $G^\psi$ are
\[
(i,\g,j)(j,\eta,k)=(i,\g\eta,k)\quad\text{and}\quad(i,\g,j)^{-1} =
(j,\g^{-1},i). 
\]
The map $\phi:G^\psi\to G$ given by $\phi(i,\g,j)=\g$ is a \hm{}
extending $\psi$.  It is easy to see that $\phi$ satisfies conditions
\partref1--\partref5 of \lemref{lem-tech} so that $\phi$
induces an equivalence of $G^\psi$
and
$G$.
\end{ex}

\begin{ex}[Pull-Back Actions]
Suppose that $\phi:H\to G$ is a groupoid \hm, and that
$\psi:=\phi\restr{\ho}$ is surjective.  Then if $\aa\in\sbrg$,
$\psi^*\A$ is an elementary \cs-bundle over $\ho$ satisfying Fell's
condition.  Furthermore, we obtain an action of $H$ on $\psi^*\A$ as
follows:
$h\cdot\(s(h),a\)=\(r(h),\alpha_{\phi(h)}(a)\)$.  We denote this
action by $\phi^*\alpha$.  Observe that
$(\psi^*\A,\phi^*\alpha)\in\sbr(H)$. 
\end{ex}

The next result show how the previous two examples are
related.

\begin{lem}
\label{lem-pullb}
Suppose that $\psi:\ho\to\go$ is a continuous open surjection, and
that $\phi:H\to G$ is a \hm{} from $H$ to $G$ and which
extends $\psi$.  If $\aa\in\sbr(G)$,
then 
\[
(\A^\grp,\alpha^\grp)\sim (\psi^*\A,\phi^*\alpha)\quad\text{in
$\sbr(H)$.}
\]
\end{lem}
\begin{proof}
The map $[u,\g,a]\mapsto \(u,\alpha_\g(a)\)$ is a $H$-equivariant
Banach bundle map.
\end{proof}

\begin{cor}
\label{cor-trivialdd}
Suppose that $\aa\in\sbrg$ and that $G$ has a Haar system
$\set{\l^u}_{u\in\go}$. Then there is a groupoid $H$ with a
Haar system
and an
$(H,G)$-equivalence $Z$ such that the Dixmier-Douady class
$\delta(\A^Z)$ is trivial in $H^3(\ho;\Z)$.
\end{cor}
\begin{proof}
Since $\A$ satisfies Fell's condition, we can find a locally finite
open cover $\U=\set{U_i}$ of $\go$ such that $\delta(\A\restr
{U_i})=0$ for all $i$.  If $X:=\coprod U_i$ and $\psi:X\to\go$ is the
usual map, then $\psi^*\A\cong\oplus\A\restr{U_i}$. 
Consequently,
$\delta(\psi^*\A)=0$.  On the other hand, we can let $H:=G^\psi$ as in
Example~\ref{ex-loctriv} so that $Z:=\grp$ is a
$(H,G)$-equivalence, and \lemref{lem-pullb}
implies that $A^\grp$ is Morita equivalent over $\ho$ to
$\psi^*(A)$ (isomorphic, in fact). Since the Dixmier-Douady
class is invariant under Morita equivalence (e.g.,
\cite[Theorem~3.5]{90c}),
$\delta(\A^\grp)=0$ as desired.

To show that $H$ has a Haar system, it will suffice to show that $Z$
has a full equivariant $s$-system (\propref{prop-haar}).  But
$Z=\coprod_i G^{U_i} =\set{(i,\g):r(\g)\in U_i}$.  Then 
\[
\alpha(f)(u):= \sum_i \int_G f(i,\g)\,d\l_u(\g)
\]
does the trick.
\end{proof}
\else

\fi
\ifdraft
%
%

\section{Generalized Twists over $H$: $\E(H)$}\label{sec-eh}

Suppose that $H$ is a locally compact groupoid.  We need to recall the
definition of a \emph{twist over $H$} or a $\T$-groupoid as defined,
for example, in \cite{alex} or \cite[\S3]{gcddi}.  Simply put, a twist
over $H$ is a principal circle bundle $j:E\to H$ equipped
with a groupoid structure such that 
$E$ is a groupoid extension of $H$ by
$\ho\times\T$:
\[
\ho\arrow{e}\ho\times\T\arrow{e,t,J}{i}E\arrow{e,t,A}{j}H.
\]
The collection of equivalence
classes of twists is a group $\Tw(H)$ with respect to Baer sum.  Thus
if $E_1$ and $E_2$ are twists over $H$, then the sum of their classes
is the class of  the quotient $E_1+E_2:= E_1\stt E_2$ of
$E_1*E_2:=\set{(e,f)\in E_1\times E_2:j_1(e)=j_2(f)}$ by the diagonal
$\T$-action
$t\cdot(e,f):=(t\cdot e,t^{-1}\cdot f)$.  

As developed in \cite{alex,renault3,gcddi}, a twist over $G$ should be
regarded as a replacement for a \emph{Borel} $2$-cocycle.  The
collection of twists forms a substitute for the groupoid cohomology
group $H^2(G,\T)$ defined in \cite{renault}.
However for our purposes, a weaker notion of equivalence of extensions
is appropriate.  Recall that strict equivalence means there is a
commutative diagram
\begin{equation}
\label{eq-equiv-*}
\begin{diagram}
\node{\ho}\arrow{e}\node{\ho\times\T}\arrow{e,t,J}{i_1}\arrow{s,=}
\node{E_1} \arrow{e,t,A}{j_1} \arrow{s,r,<>}{\phi} \node{H}
\arrow{s,=} \\
\node{\ho}\arrow{e}\node{\ho\times\T}\arrow{e,t,J}{i_2}\node{E_2}
\arrow{e,t,A}{j_2} \node{H}
\end{diagram}
\end{equation}
where $\phi$ is a groupoid isomorphism.  That is, $E_1$ and $E_2$ are
isomorphic over $H$.

 We want to weaken this notion in two steps.  In
rough terms, the first is to replace the isomorphism $\phi$ in
\eqref{eq-equiv-*} with an equivalence of $\T$-groupoids while
maintaining equality in the other vertical positions.  The second step
will be to allow the groupoid $H$ to vary over equivalent groupoids. 
To make these statements precise will take a bit of work.  First
recall that two twists $E$ and $F$ over $H$ are equivalent as
$\T$-groupoids if there is an $(E,F)$-\teq{}
as defined in \cite[Definition~3.1]{gcddi}; that is, if there is
an ordinary $(E,F)$-equivalence $Z$ for which the $\T$-actions on $Z$
induced by $E$ and $F$ coincide.  It is straightforward
to check that the quotient $Z/\T$ is naturally an $(H,H)$-equivalence.
We 
define two twists $E$ and $F$ to be \emph{equivalent over $H$} if
there is a \efteq{} $Z$ such that $Z/\T$ is isomorphic to $H$ (as
$(H,H)$-equivalences).  In this case we write $E\simh F$, and say
that $E$ and $F$ are \emph{equivalent over $H$}.

If $[E]=[F]$ in $\Tw(H)$, then there is a \emph{twist isomorphism}
\[
\begin{diagram}
\node{E}\arrow[2]{e,t}{\phi} \arrow{se,b}{j_1} \node[2]{F}
\arrow{sw,b}{j_2} \\
\node[2]{H}
\end{diagram}
\]
which is a bundle isomorphism as well as a groupoid \hm.
In particular, $\grp$ is an \efteq, and $\grp/\T$ satisfies
$j_1(e)[s(e),f] = [r(e),\phi(e)f]$.  Since $j_1(e)=j_2\(\phi(e)\)$,
the map $[u,f]\mapsto \(u,j_2(f)\)$ is an isomorphism of
$\grp/\T$ onto
$\Gr(\id_H)$.  Thus $E\simh F$ in this case.  The converse is not
valid (see Remark~\ref{rem-7.4}).   To examine this, consider the group
\hm{}
$\epsilon:H^1(\ho,\sT)\to
\Tw(H)$ defined on page~124 of \cite{gcddi}.  Recall that $\epsilon$
is defined as follows.  We identify $H^1(\ho,\sT)$ with isomorphism
classes of principal $\T$-bundles over $\ho$.  If $p:\tb\to\ho$ is
such a bundle, we let $\btb$ be the conjugate bundle.  Note that
viewing $\tb$ as a right $\T$-space, $\btb$ is the opposite left
$\T$-space: $t\cdot\iota(\xi)=\iota(\xi\cdot t^{-1})$.  Then let
$\tbstbtb$ be the $\T$-bundle over $\ho\times\ho$ obtained as the
quotient of $\tb\times\btb$ by the diagonal $\T$-action:
$t\cdot\(\xi,\iota(\eta)\) = \(\xi\cdot t,t^{-1}\cdot\iota(\eta)\)$.
Then $\epsilon\([\tb]\)$ is the class of the pull-back
$\pi^*(\tbstbtb)$ where $\pi:H\to\ho\times\ho$ is the map
$\pi(h)=\(r(h),s(h)\)$. 

If $E$ is a principal $\T$-bundle over $H$ and $\tb$ is a principal
$\T$-bundle over 
$\ho$, then the quotient $\tb\stt E$ of
$\tb*E:=\set{(\xi,e)\in\tb\times E: p(\xi)=r(e)}$ by the $\T$-action
$t\cdot (\xi,e):=(\xi\cdot t,t^{-1} e)$ is also a principal
$\T$-bundle over $H$.  There is a similar quotient $E\stt\tb$ of
$E*\tb:=\set{ (e,\eta)\in E\times\tb:s(e)=p(\eta)}$.  Our interest in
these constructions is that if $E$ is also a twist over $H$, then the
same is true of 
\begin{equation}
\label{eq-tbtwist}
\tb\stt E\stt\btb:=(\tb\stt E)\stt \btb\cong \tb\stt(E\stt\btb);
\end{equation}
we have $r\([\xi,e,\iota(\eta)]\)=p(\xi)$, 
$s\([\xi,e,\iota(\eta)]\)=p(\eta)$, and
$[\xi,e,\iota(\eta)][\eta\cdot t,f,\iota(\zeta)]=[\xi,tef,\iota(\zeta)]$.  

\begin{lem}
\label{lem-silly}
Suppose that $E$ is a twist over $H$ and that $\tb$ is a principal
$\T$-bundle over $\ho$.  Then
\[
[E]+\epsilon\([\tb]\)=[\tb\stt E\stt\btb].
\]
\end{lem}
\begin{proof}
It is not hard to see that $[E]+\epsilon\([\tb]\)$ is represented by
the quotient of $\set{(e,h,\xi,\iota(\eta))\in E\times
H\times\tb\times\btb:
\text{$j(e)=h$, $r(h)=p(\xi)$, and $s(h)=p(\eta)$}}$ by the
$\T^2$-action $(t,s)\cdot(e,h,\xi,\iota(\eta))= (e\cdot
t^{-1},h,\xi\cdot ts,s^{-1}\iota(\eta))$.  Then $(e,h,\xi,\iota(\eta))
\mapsto (\xi,e,\iota(\eta))$ induces a twist isomorphism.
\end{proof}

\begin{lem}
\label{lem-alex}
Suppose that $E$ and $F$ are $\T$-groupoids over $H$.  Then $E\simh F$
if and only if there is a principal $\T$-bundle over $\ho$ such that 
\[
[F]=[\tb\stt E\stt \btb].
\]
\end{lem}
\begin{proof}
Suppose that $E\simh F$ and that $Z$ is an \efteq{} such that $Z/\T$
is isomorphic to $H$ as an $(H,H)$-equivalence.  Composing with the
orbit map then gives $Z$ the structure of a $\T$-bundle over $H$,
$p:Z\to H$ and the restriction to $\ho$ is therefore a principal
$\T$-bundle $\tb$ over $\ho$.  Since $Z/\T\cong H$, we have
$s\(p(z)\)=s(z)$ and $r\(p(z)\)=r(z)$.  Therefore we obtain a map from
$\tb*E=\set{(\xi,e)\in\tb\times E:p(\xi)=r(e)}$ to $E$ by restricting
the action map from $Z*E\to Z$.  The restriction is invariant with
respect to the given $\T$-action on $\tb*E$ and defines a continuous
map $\phi$ from $\tb\stt E$ to $Z$.  It is not hard to see that $\phi$
is a bundle map, and hence a bundle isomorphism (cf., e.g.,
\cite[Theorem~3.2]{husemoller}).  Therefore $F$ is groupoid isomorphic
to the imprimitivity groupoid
$E^{\tb\stt E}=(\tb\stt E)\ste
\opposite{(\tb\stt E)}$. 
The latter is clearly the quotient of $\set{(\xi,e,\eta,f)\in
\tb\times E\times \tb\times E:\text{$p(\xi)=r(e)$, $p(\eta)=r(f)$, and
$s(e)=s(f)$}}$ by the action of $\T^2\times E$: $(t,s,g)\cdot
(\xi,e,\eta,f) = (\xi\cdot t, t^{-1}\cdot eg,\eta\cdot s,s^{-1}\cdot
fg)$.  The map $[\xi,e,\eta,f]\mapsto [\eta,ef^{-1},\iota(\eta)]$ is a
groupoid isomorphism onto $\tb\stt E\stt\btb$.
On the other hand, given $f\in F$ and $\xi\in\tb$ with $r(\xi)=s(f)$,
the isomorphism of $F$ with $Z\ste\opposite{Z}$ takes $f$ to
$[f\cdot\xi,\xi]$.  Choose $\eta\in\tb$ and $e_f\in E$ such that
\begin{equation}
\label{eq-eqjs}
\eta\cdot e_f=f\cdot\xi.
\end{equation}
  Then the induced isomorphism of $F$ with
$\tb\stt E\stt\btb$ takes $f$ to $[\eta,e_f,\xi]$.  But
\eqref{eq-eqjs} and the fact that $Z/\T\cong H$ implies that
\[
j_F(f)=j_F(f)\cdotp(\xi)=p(\eta)\cdot j_E(e_f)=j_E(e_f).
\]
It follows that the isomorphism of $F$ with $\tb\stt E\stt\btb$ is
also a bundle isomorphism.  This proves ``only if'' implication.

On the other hand, if $[F] = [\tb\stt E\stt\btb]$ for some $\T$-bundle
over $\ho$, then $F$ is the imprimitivity groupoid for $\tb\stt E$.
Then $(\tb\stt E)/\T\cong \Gr(\id_H)$.  This completes the proof.
\end{proof}

We will denote the trivial twist over $H$ by $H\times\T$.

\begin{cor}
\label{cor-simh}
The collection of classes $[E]$ in $\Tw(H)$ such that $E\simh
H\times \T$ is a subgroup $W$ of $\Tw(H)$.  In particular,
$\simh$ is an equivalence relation on $\Tw(H)$, and
$\E(H):=\Tw(H)/\simh$ can be identified with the quotient group
$\Tw(H)/W$.  Addition in $\E(H)$ is defined by
\[
\lbb E \rbb + \lbb F \rbb = \lbb E\stt F \rbb.
\]
\end{cor}
\begin{proof}
By the previous lemma, $W=\epsilon\(H^1(\ho,\sT)\)$, and
$\simh$ equivalence classes coincide with $W$-cosets.
\end{proof}

\begin{remark}
\label{rem-7.4}
The subgroup $W$ in \corref{cor-simh}
can be nontrivial.  For example, if $\psi:X\to T$ is a local
homemorphism and if $H$ is the corresponding equivalence relation
$R(\psi)$ in $X\times X$, then the six term exact sequence in either
\cite{alex2} or \cite{gcddi} implies that $W$ is the cokernal of the
pull-back $\psi^*:H^1(T,\sT)\to H^1(X,\sT)$.  This can certainly be
nontrivial; for example, let $X=T=\T^2$ and $\psi$ a double cover.
\end{remark}

\else

\fi
\ifdraft
%
%

\section{The Second Isomorphism Result}\label{sec-e(x)}

The goal of this section is \thmref{thm-main2}, which identifies
$\E(G)$ with the special subgroup of $\brg$ consisting of those
elementary $G$-bundles with vanishing Dixmier-Douady invariant.

It is well-known that Morita equivalence classes of continuous-trace
\cs-algebras with spectrum $T$ are parameterized by classes in
$H^3(T;\Z)$ (cf., e.g., \cite[Theorem~3.5]{90c}).  If $A=\gsuboT(\A)$
is a continuous-trace \cs-algebra with spectrum $T$, then the
corresponding class in $H^3(T;\Z)$ is denoted by $\delta(A)$ (or
$\delta(\A)$) and is called the
\emph{Dixmier-Douady class of $A$ (or $\A$)}.  In \cite{90c},
$\delta(A)$ is constructed as the obstruction to $A$ being Morita
equivalent to $C_0(T)$.  The later is equivalent to the existence of a
Hilbert $C_0(T)$-module $\H$ such that $A\cong\K(\H)$.  The classical
interpretation of $\delta(\A)$ is as the obstruction to the bundle
$\A$ arising as the bundle $\kkhh$ associated to a Hilbert bundle
$\HH$ (\cite[Theorem~10.7.15]{dix}).  Of course, if $\H=\gsuboT(\HH)$
and $A=\gsuboT(\A)$, then $\kh=\gsuboT(\kkhh)$.  Despite the
fourth author's fondness for the methods of \cite{90c}, the bundle
interpretation is the more natural one here.  However, as the
Dixmier-Douady class is invariant under Morita equivalence (over the
spectrum), we note that if $\A$ and $\B$ are Morita equivalent
\cs-bundles satisfying Fell's condition, then
$\delta(\A)=\delta(\B)$.  This allows us to make the following
definition.

\begin{definition}
If $G$ is a \sclcg, then
\[
\sbrzg:=\set{\aa\in\sbrg:\delta(\A)=0},
\]
and
\[
\brzg:=\set{[\A,\a]\in\brg:\aa\in\sbrzg}.
\]
\end{definition}
$\brzg$ is a subgroup in view of \cite[Proposition~3.3]{ckrw}.

If $\hs$ is a Hilbert space, we will always view $\Aut\(\K(\hs)\)$ as a
topological group with the point-norm topology.  Of course, the map
$U\mapsto \Ad(U)$ identifies $\Aut\(\K(\hs)\)$ with the quotient of the
unitary group $U(\hs)$, equipped with the strong operator topology by
its center $\T\cong \T\cdot1_{\hs}$.  We will use $\eltwo$ as the
generic infinite dimensional separable Hilbert space and write $\K$
for $\K(\eltwo)$.  Note that $\KK(\go\times\eltwo)$ is the trivial
bundle $\go\times\K$.

We will write $\rg$ for the set of continuous groupoid \hm s $\pi:G
\to \Aut(\K)$.  As observed in \cite[\S4]{alex3}, 
each element $\pi\in\rg$ defines an element of $\sbrzg$
\begin{equation}
\label{eq-bar*}
(\go\times\K,\bar\pi)
\end{equation}
in $\sbrzg$,
where
\[
\bar\pi_\g\(s(\g),T\):=\(r(\g),\pi_\g(T)\).
\]
Clearly, if $(\go\times\K,\alpha)\in\sbrzg$, then $\alpha=\bar\rho$
for some $\rho\in\rg$.
It is a pleasant surprise that all elements of $\brzg$ have a
representative of the form of \eqref{eq-bar*}.

\begin{lem}
\label{lem-pi}
If $G$ is a \sclcg{} and if $\aa\in\sbrzg$, then there is a
$\pi\in\rg$ such that
\[
\aa\sim(\go\times\K,\bar\pi).
\]
\end{lem}
\begin{proof}
Since $\delta(\A)=0$, there is a Hilbert bundle $\HH$ such that
$\A=\kkhh$ (\cite[Theorem~10.7.15]{dix}).  Let $\H=\gsubogo(\HH)$.
Since $A=\gsubogo(\A)$ is separable, $\H$ is countably generated and
\cite[Proposition~7.4(ii)]{lance} implies that the Hilbert
$C_0(\go)$-modules $\H\tensor\eltwo$ and $\H_{C_0(\go)}:=
C_0(\go)\tensor \eltwo$ are isomorphic.  It follows from
Remark~\ref{rem-**} that 
\begin{equation}
\label{eq-keyiso}
\HH\tensor(\go\times\eltwo)\cong \go\times\eltwo
\end{equation}
as Hilbert bundles over $\go$.  Since $\tau_{\go}\sim
(\go\tensor\K,\tau\tensor 1)$ acts as the identity in $\brg$ and
since
\[
\kkhh\tensor(\go\times\K)\cong \KK\(\HH\tensor(\go\times\eltwo)\),
\]
if follows that 
\begin{align*}
\aa &= \(\kkhh,\alpha\) \sim \(\kkhh,\alpha\) \\
&\sim \(\KK(\HH\tensor(\go\times\eltwo)\),\alpha\tensor(\tau\tensor1)\),
\end{align*}
which is covariantly isomorphic to
$\(\KK(\go\times\eltwo),\tilde\alpha\)$ for some action
$\tilde\alpha$.  As remarked above, $\tilde\alpha$ must be of the form
$\bar\pi$ for some $\pi\in\rg$.
\end{proof}

If $\pi\in\rg$, then there is an associated twist over $G$ defined by
\[
E(\pi):=\set{(\g,U)\in G\times U(\eltwo): \pi_\g=\Ad(U)}.
\]
The main result of this section is the following.

\begin{thm}
\label{thm-main2}
Suppose that $G$ is a \sclcg.  Then
\[
[\go\times\K,\bar\pi]\mapsto \eclass{E(\pi)}
\]
defines an isomorphism $\thetag$ of $\brzg$ onto $\E(G)$.
\end{thm}

The proof is fairly involved and is broken into a number of
intermediate results.  Naturally, the main tool is the relationship
between Morita equivalence of systems and the elements of $\rg$.  For
this we need a couple of definitions.

We say that $\pi$ and $\rho$ in $\rg$ are \emph{exterior equivalent}
if there is a continuous map $u:G\to U(\eltwo)$ such that
\begin{enumerate}
\item
$\rho_\g=\Ad(u_\g)\circ \pi_\g$ for all $\g\in G$, and
\item
$u_{\g\eta}=u_\g\pi_\g(u_\eta)$ for all $(\g,\eta)\in G^{(2)}$.
\end{enumerate}
In this case, we write $\pi\simext\rho$ and call $u$ a $1$-cocycle
implementing this equivalence.

It is not hard to see that $E(\pi)$ and $E(\eta)$ are isomorphic
elements of $\Tw(G)$ if and only if $\pi\simext\rho$.  To see this,
suppose that $u$ is a $1$-cocycle as above.  Then $(\g,V)\mapsto
(\g,u_\g V)$ is a continuous twist isomorphism of $E(\pi) $ onto
$E(\rho)$.  Conversely, if $\phi:E(\pi)\to E(\rho)$ is a twist
isomorphism, then since $U(\eltwo)$ is a group, 
\[
\phi(\g,V)=(\g,u_\g V)
\]
for some continuous function $u:G\to U(\eltwo)$.  Since $\phi$ is also
a groupoid \hm, it is immediate that
$u$ is a $1$-cocycle implementing an exterior equivalence between
$\pi$ and $\rho$.

We say that $\pi$ and $\eta$ are \emph{cohomologous} (written
$\pi\simcoh \rho$) if there is a continuous function
$\phi:\go\to\Aut(\K)$ such that
\[
\pi_\g\circ\phi_{s(\g)} = \phi_{r(\g)}\circ \rho_\g.
\]
It is not hard to see that $(\go\times\K,\bar\pi)$ and
$(\go\times\K,\bar\rho)$ are covariantly isomorphic if and only if
$\pi\simcoh \rho$.

The main tool in the proof of \thmref{thm-main2} is the following.

\begin{prop}
[{\cite{masuda}}, {\cite[Proposition~14]{alex3}}]%
\label{prop-alex14}
Suppose that $G$ is a \sclcg, and that $\pi,\rho\in\rg$.  Then
\[
(\go\times\K,\bar\pi)\sim (\go\times\K,\bar\rho)
\]
if and only if there is a $\sigma\in\rg$ such that 
\[
\pi\simext\sigma\quad\text{and}\quad\sigma\simcoh\rho.
\]
\end{prop}
\begin{proof}
The ``if'' direction is clear.  For the ``only if'' direction, let
$\xv$ implement an equivalence between $(\go\times\K,\bar\pi)$ and
$(\go\times\K,\bar\rho)$, and let $\X=\gotoco(\go,\bX)$,
 $\A=\go\times\K$, and
$A=\gsubogo(\A)=C_0(\go,\K)$.  Then $\X$ is an $A\smeover\go A$-\ib.
By \cite[Proposition~3.4]{bgr}, $\X$ is isomorphic (as an $A\sme
A$-\ib) to $A_\kappa$ where $\kappa\in\Aut(A)$.  Since $\X$ is an \ib{}
over $\go$, the same must be true of $A_\kappa$ (using
\cite[Proposition~1.11]{ra1}, for example).  Then
$\kappa\in\Aut_{C_0(\go)}(A)$ and must be given by a continuous
function $\phi:\go\to\Aut(\K)$.  Therefore we may as well assume that
$\bX=\go\times\K$ and that 
\begin{gather*}
\brip\A<{(u,T)},(u,S)>=\(u,\phi_u(T^*S)\)\qquad
\blip\A<{(u,T)}, (u,S)> = (u,TS^*) \\
\brip\A<{\vg\(s(\g),T\)}, \vg\(s(\g),S\)> = \bar\rho_\g\( \brip\A <
{\(r(\g),T\)} , \(r(\g),S\)> =\(r(\g), \rho_\g\(\phi_{r(\g)}(T^*S)\)\)
\\
\blip\A<\vg\({s(\g),T}\),\vg\(s(\g),S\)> = \bar\pi_\g\( \blip\A
<{\(r(\g), T\)},\(r(\g),S\)> = \(r(\g),\pi_\g(TS^*)\),
\end{gather*}
with appropriate left and right actions.

Recall that we have an isomorphism $V:s^*A\to r^*A$ defined by
$V(f)(\g):= \vg\(f(\g)\)$.  Then it is straightforward to compute that
$V(fg)=\tilde\pi(f)V(g)$, where $\tilde\pi(f)(\g)=\bar\pi_\g\(f(\g)\)$.
Thus we can define $L:r^*A\to r^*A$ by $L(f):=
V\(\tilde\pi^{-1}(f)\)$.  Note that by construction, $L(fg)=fL(g)$.
Of course, $r^*A$ is isomorphic to $C_0(G,\K)$ and is a left Hilbert
$C_0(G,\K)$-module.  Furthermore, $L$ is an adjointable operator on
$r^*A$:
$L^*(f)(\g):= \bar\pi_\g^{-1}\(V_{\g^{-1}}\(f(\g)\)\)$.  However
$\L(r^*A)$ is isomorphic to $M\(\K(r^*A)\)\cong M\(C_0(G,\K)\)$ by
\cite[Lemma~16]{green1}\footnote{See also Theorem~2.4 of \cite{lance},
where the result is attributed to Kasparov.}.  But
$M\(C_0(G,\K)\)\cong C_b\(G,B(\eltwo)_{*-s}\)$, where
$B(\eltwo)_{*-s}$ denotes the bounded operators on $\eltwo$ with the
$*$-strong topology \cite[Corollary~3.4]{apt}\footnote{Formally, we
should have the strict topology on $B(\eltwo)$, but these topologies
coincide on \emph{bounded} subsets.}.  Therefore $L(f)=f\cdot u^*$ for
some $u\in C_b\(G,B(\eltwo)_{*-s}\)$, where $f\cdot
u^*(\g)=\(r(\g),Tu_\g^*\)$ if $f(\g)=\(r(\g),T\)$.  Alternatively, we
can define $\tilde U^*\in r^*A$ by $\tilde U_\g=\(r(\g),u_\g^*\)$ and
$L(f)=f\tilde U^*$.  In any case, we can compute that for all $f\in
r^*A$
\[
\blip r^*A<L(f),L(f)>=\lip r^*A<f,f>= \blip r^*A<L^*(f),L^*(f)>,
\]
while on the other hand
\[
\blip r^*A<L(f),L(f)>(\g)=f(\g)\tilde U_\g \tilde U^*_\g f(\g)^*.
\]
Therefore $u_\g u_\g^*=I$.  Similarly, $u_\g^*u_\g=I$ and $u$ is a
continuous function from $G$ into $U(\eltwo)$.  Since
$V(f)=\tilde\pi(f)\tilde U^*$ and since $V$ is an action, it follows
that $u_{\g\eta}=u_\g\pi_\g(u_\eta)$.  Thus if we define
$\sigma\in\rg$ by
\[
\sigma_\g:=\Ad(u_\g)\circ \pi_\g.
\]
then $\pi\simcoh\sigma$ by construction.  

If $f,g\in\X=\gsubogo(\bX)$, then
\begin{align*}
\brip\A<{\vg\(f(\g)\)},\vg\(g(\g)\)> &= \bar\phi_{r(\g)}
\(\vg\(f\(s(\g)\)\)^* \vg\(g\(s(\g)\)\)\) \\
&= \bar\phi_{r(\g)}\(\tilde U_\g \bar\pi_\g \(f\(s(\g)\)^*g\(s(\g)\)
\) \tilde U^*_\g\) \\
&
= \bar\phi_{r(\g)} \(\bar\sigma_\g\( f\(s(\g)\)^*g\(s(\g)\) \) \).
\end{align*}
Since the left-hand side is also equal to
$\bar\rho_\g\circ\bar\phi_{s(\g)} \(f\(s(\g)\)^*g\(s(\g)\)\)$, it
follows that $\phi_{r(\g)} \circ \sigma_\g = \rho_\g \circ
\phi_{s(\g)}$.  That is, $\sigma\simcoh\rho$; this completes the proof.
\end{proof}

\begin{cor}
\label{cor-extcoh}
If $G$ is a \sclcg, if $\pi,\rho\in\rg$, and if
\[
(\go\times\K,\bar\pi)\sim(\go\times\K,\bar\rho),
\]
then $\eclass{E(\pi)}=\eclass{E(\rho)}$ in $\E(G)$.
\end{cor}

\begin{proof}
Using the proposition, there is a $\sigma\in\rg$ such that
$\pi\simext\sigma$ and $\sigma\simcoh\rho$.  Since exterior
equivalence implies that $E(\pi)$ and $E(\sigma)$ are isomorphic as
twists, it will suffice to show that $E(\sigma)\simh E(\rho)$.    Let
$\phi:\go\to \Aut(\K)$ be such that 
\[
\sigma_\g\circ\phi_{s(\g)} = \phi_{r(\g)}\circ \rho_\g.
\]
Let 
\[
\Lambda:=\set{(u,W)\in \go\times U(\eltwo):\phi_u=\Ad(W)}.
\]
Then $\Lambda$ is a principal circle bundle over $\go$, and
$\epsilon(\Lambda)$ is the class of 
\[
\set{[U,\g,V]\in U(\eltwo)\times G\times U(\eltwo)/\T:
\text{$\Ad(U)=\phi_{r(\g)}$ and $Ad(V)=\phi_{s(\g)}$}},
\]
where $t\cdot(U,\g,V)=(tU,\g,tV)$.  Then
\[
\bigl[[U,\g,V],(\g,W)\bigr]\mapsto (\g,UWV^*)
\]
is a twist isomorphism of $\epsilon(\Lambda)+E(\rho)$ onto
$E(\sigma)$.  This completes the proof (\corref{cor-simh}).
\end{proof}

\begin{cor}
\label{cor-injhm}
The map $\thetag\([\go\times\K,\bar\pi]\)=\eclass{E(\pi)}$ is a
well-defined injective \hm{} from $\brzg$ to $\E(G)$.
\end{cor}

\begin{proof}
\lemref{lem-pi} imples every class in $\brzg$ has the form
$[\go\times\K,\bar\pi]$, and $\thetag$ is well-defined by
\propref{prop-alex14}.  To see that $\thetag$ is a \hm, we need only
consider
\begin{equation}
\label{eq-prod*}
(\go\times\K,\bar\pi)\tensor(\go\times\K,\bar\rho) = (\go\times
\K\tensor\K,\bar \pi\tensor\bar\rho).
\end{equation}
Fix a Hilbert space isomorphism $W:\eltwo\tensor\eltwo\to\eltwo$, and
let $\phi:=\Ad(W):\K\tensor\K\to \K$ be the corresponding isomorphism.
Then \eqref{eq-prod*} is covariantly isomorphic to
$(\go\times\K,\bar\sigma)$ where $\sigma_\g=\phi\circ (\pi_\g\tensor
\rho_\g)\circ\phi^{-1}$.  Then
\[
E(\pi)+E(\rho)= \set{[U,\g,V]\in U(\eltwo)\times G\times U(\eltwo)/\T:
\text{$\Ad(U)=\pi_\g$ and  $\Ad(V)=\rho_\g$}},
\]
where $t\cdot(U,\g,V)=(tU,\g,t^{-1}V)$.  Therefore
\[
[U,\g,V]\mapsto (\g,W(U\tensor B)W^*)
\]
is a twist isomorphism of $E(\pi)+E(\rho)$ onto $E(\sigma)$.  It
follows that $\thetag$ is an \hm.

Thus to see that $\thetag$ is injective, it suffices to show it has
trivial kernel.  However if $\eclass{E(\pi)}=0$, then
$E(\pi)=\epsilon(c)$ for some $c\in H^1(\go,\sT)$.  But by
\cite[Theorem~2.1]{pr1}, there is a continuous function
$\sigma:\go\to\Aut(\K)$ such that
$\a\in\Aut_{C_0(\go)}\(C_0(\go,\K)\)$ defined by
$\alpha(f)(u)=\sigma_u\(f(u)\)$ has Phillips-Raeburn obstruction
$\zeta(\alpha) = c$.  Then $c=[\Lambda]$ where
\[
\Lambda=\set{(u,W)\in\go\times U(\eltwo):\Ad(U)=\sigma_u}.
\]
Then 
\[
\epsilon(\Lambda)=\set{[R,\g,S]\in U(\eltwo)\times G\times U(\eltwo)/\T:
\text{$\Ad(R)=\sigma_{r(\g)}$ and $\Ad(S)=\sigma_{s(\g)}$}},
\]
where $t\cdot(R,\g,S)=(tR,\g,tS)$.
If we define $\rho\in\rg$ by
$\rho_\g=\sigma_{r(\g)}\sigma_{s(\g)}^{-1}$, then
\[
[R,\g,S]\mapsto (\g,RS^*)
\]
is a twist isomorphism of $E(\pi)=\epsilon(\Lambda)$ onto $E(\phi)$.
Since $\rho$ is cohomologous to the trivial \hm, the result follows
from (the easy half) of
\propref{prop-alex14}. 
\end{proof}

This leaves only surjectivity.  At this point, our methods require the
existence of a Haar system.

\begin{prop}
\label{prop-surjectivity}
Suppose that $G$ is a \sclcg{} with a Haar system
$\set{\l^u}_{u\in\go}$.  Then if $E\in\Tw(G)$, there is a $\pi\in\rg$
such that $E\cong E(\pi)$.
\end{prop}

Recall from \cite{alex3} 
that a \emph{twist representation} of $E\in\Tw(G)$ is a pair
$(\HH,U)$ consisting of a Hilbert bundle $\HH$ and an action $U$ of
$E$ on $\HH$ such that
\[
U_{t\cdot e}=tU_e \quad\text{for all $t\in \T$}.
\]
Then $\alpha_{j(e)}:= \Ad(U_e)$ is a well-defined action of $G$ on
$\kkhh$, and $\(\kkhh,\alpha\)\in\sbrzg$.  Naturally we say that
$(\HH,U)$ is a \emph{separable} twist representation if
$\H=\gsubogo(\HH)$ is separable.

\begin{lem}
\label{lem-one}
Suppose that $G$ is a \sclcg{} and that $E\in\Tw(G)$.  If $(\HH,U)$ is
a separable twist representation of $E$ and if $\(\kkhh,\alpha\)$ is the
associated element of $\sbrzg$, then
\[
\thetag\([\kkhh,\alpha]\)=\eclass E.
\]
\end{lem}

\begin{proof}
It follows from \eqref{eq-keyiso} that there is a bundle isomorphism
$V=\set{V_u}$ of $\HH\tensor(\go\times\eltwo)$ onto $\go\times\eltwo$
which intertwines $U\tensor(\tau\tensor1)$ with an action $\tilde U$
on the trivial bundle.  Note that 
\[
\tilde U_e=V_{r(e)}\(U_e\tensor(\tau_e\tensor 1)\)V_{s(e)}^*.
\]
Thus $\tilde U$ is continuous and there exists a continuous \hm{}
$W:E\to U(\eltwo)$ such that 
\[
\tilde U_e\(s(e),\xi\)=\(r(e), W_e(\xi)\).
\]
This forces $W_{t\cdot e}=tW_e$.  Then $\pi_{j(e)}:= \Ad(W_e)$ is an
element of $\rg$ and 
\[
\(\kkhh,\alpha\)\sim (\go\times\K,\bar\pi).
\]
It is immediate that $e\mapsto \(j(e),W_e\)$ is a twist isomorphism of
$E$ and $E(\pi)$.
\end{proof}

\begin{proof}
[Proof of \propref{prop-surjectivity}]
At this point, all we need to do is to produce a twist representation
of $E\in\Tw(G)$.
As in \cite{mw3}, we consider 
\[
C_c(G;E):=\set{f\in C_c(E): f(t\cdot e)=tf(e)}.
\]
Notice that $C_c(G;E)$ is a pre-Hilbert $C_0(\go)$-module:
\begin{equation}
\label{eq-dagger}
\rip C_0(\go)<f,g>(u):= \int_G \overline{f(e)}g(e)\,d\l^u\(j(e)\).
\end{equation}
The completion is a Hilbert $C_0(\go)$-module $\H=\gsubogo(\HH)$,
where $\HH$ is a Hilbert bundle with fibres $\H_u$ equal to the
completion of $C_c\(r^{-1}(u)\)$ with respect to the restriction of
\eqref{eq-dagger}.  If $e\in E$, then define
$U_e:C_c\(r^{-1}\(s(e)\))\to C_c\(r^{-1}\(r(e)\)\)$ by
$U_e(f)(e')=f(e^{-1}e')$.  Then $U_e$ is easily seen to extend to a
unitary isomorphism of $\H_{s(e)}$ onto $\H_{r(e)}$ such that
$U_{ee'}=U_e\circ U_{e'}$.  To see that $U$ is a continuous action of
$E$, we want to apply \lemref{lem-1.12} to $U:s^*\H=\gotoco(G,s^*\HH)
\to r^*\H=\gotoco(G,r^*\HH)$.
Note that $\Gamma:=C_c(G;E)$ can be viewed as a dense subspace of
sections in $\H$.  If $f\in C_c(G;E)$ and $\phi\in C_c(E)$, then
\[
\phi\tensor s^*f(e):= \(e,\phi(e)f\(s(e)\)\)
\] 
is in $s^*\H$ and the collection $s^*\Gamma$ of such sections is dense
in $s^*\H$ by \cite[Proposition~II.14.6]{fell-doran}.

Let $C_c(G;E^{(2)})$ be the pre-Hilbert $C_0(E)$-module consisting of
continuous compactly supported functions on $E^{(2)}=\set{(e,e')\in
E\times E: s(e)=r(e')}$ such that $f(e,t\cdot e')=tf(e,e')$ with 
\[
\rip C_0(E)<f,g>(e) := \int_G \overline{f(e,e')} f(e,e') \,d\l^{s(e)}(e').
\]
Each element of $s^*\Gamma$ defines an element of $C_c(G;E^{(2)})$:
$\phi\tensor s^*f(e,e')=\phi(e)f(e')$.  The image of $s^*\Gamma$ is
therefore dense in $C_c(G;E^{(2)})$ by
\cite[Proposition~II.14.6]{fell-doran}.  This allows us to identify
$s^*\H$ with the completion of $C_c(G;E^{(2)})$.  Similarly, we can
identify $r^*\H$ with the completion of $C_c(G;E\understar s E)$ where
$E\understar sE:=\set{(e,e')\in E\times E: s(e)=s(e')}$.  Therefore
$L(f)(e,e'):= f(e,e^{-1}e')$ extends to an isomorphism of $s^*\H$ onto
$r^*\H$.  It follows that $U$ is a continuous action as claimed.
It only remains to observe that $U_{t\cdot e}=tU_e$; thus, $(\HH,U)$
is a twist representation.
\end{proof}

\begin{proof}
[Proof of \thmref{thm-main2}]
$\thetag$ is an injective \hm{} by \corref{cor-injhm} and surjective
by \propref{prop-surjectivity}.
\end{proof}

In \secref{sec-last}, we will show that $\brg$ can be computed in
terms of the $\E(H)$ for groupoids $H$ equivalent to $G$.  However if
$\go$ is a reasonable space, it is easy to recover $\brg$ from the
generalized twists over a single equivalent groupoid.

\begin{cor}
\label{cor-singlecover}
Suppose that $G$ is a \sclcg{} with a Haar system, and that $\go$ is
locally contractible.  Then there is an equivalent groupoid $H$ such
that $\brg$ is isomorphic to $\E(H)$.
\end{cor}

\begin{proof}
Let $\U=\set{U_i}_{i\in I}$ be a locally finite cover of $\go$ by
contractible open sets.  In fact, any cover for which
$H^3(U_i;\Z)=\set0$ for all $i$ will do.  Anyway, if $\aa\in\sbrg$,
then $\delta(\A\restr{U_i})=0$ for all $i$.  Then if $X:=\coprod U_i$,
$\psi:X\to\go$, and $\phi:G^\psi\to G$ are as in
\corref{cor-trivialdd}, then $Z:=\grp$ is a \equi(G^\psi,G) and
$\phi^Z\([\A,\a]\) \in \brz(G^\psi)$.  Since $\phi^Z$ is an
isomorphism, $\brz(G^\psi)=\br(G^\psi)$.  Now the result follows from
\thmref{thm-main2}. 
\end{proof}
\else

\fi
\ifdraft
%
%

\section{$\extgt$}\label{sec-four}

As we saw in \secref{sec-eh}, 
$\simh$ is indeed a weaker notion of equivalence of extensions
and is the ``first step'' alluded to earlier.  The next step is to
allow the groupoid $H$ to vary over equivalent groupoids.  The idea is
that extensions of $H'$ and $H''$ should be compared in $\E(H)$ where
$H$ is a ``larger'' groupoid equivalent to $H'$ and $H''$.  This
process is formalized by taking an inductive limit.  It turns out to
be appropriate to work with principal $G$-spaces rather than
equivalent groupoids.

Let $\pg$ be the collection of second countable\footnote{We will
eventually need separability, and therefore second countability, for
some of our results.  At the moment, some restriction on the ``size''
of the elements of $\pg$ is necessary in order that it be a set.}
principal right $G$-spaces which have a full equivariant $s$-system.  
By \propref{prop-haar}, these are precisely the principal $G$-spaces whose
imprimitivity groupoid has a Haar system.

Motivated by
Example~\ref{ex-4.1}, if $X,Y\in\pg$, then we will say that \emph{$X$
is a refinement of $Y$} if there is a continuous, open,
$G$-equivariant surjection $\phi:X\to Y$.  In this event we will write
$Y\preceq X$, or $Y\preceq_\phi X$ if we want to keep the map $\phi$ in
evidence.  Note that $\pg$ is a directed set with respect to
$\preceq$: given $X,Y\in\pg$, let $Z:=X*Y=\set{(x,y)\in X\times Y:
s(x)=s(y)}$.  Then $Z$ is a principal $G$-space with respect to the
diagonal action.  If $\alpha$ and $\beta$ are full equivariant
$s$-systems on $X$ and $Y$, respectively, then
\[
\pi(f)(u):=\int_X\int_Y f(x,y)\,d\beta_u(y)\,d\alpha_u(x)
\]
is a full equivariant $s$-system for $Z$; that is, $Z\in\pg$.
Using the projection maps, we see
that that $X\preceq Z$ and $Y\preceq Z$.

If $Y\preceq_\phi X$, then there is an associated groupoid \hm{}
$\tphi:G^X\to G^Y$ given by
$\tphi\(\eq[x,z]\)=\eq[\phi(x),\phi(z)]$.  Using
\lemref{lem-tech}, it is not hard to see that $\tphi$ induces an
equivalence between $G^X$ and $G^Y$.\footnote{In fact, the equivalence
is isomorphic to $X\stg\opposite{Y}$ which is independent of $\phi$
(\lemref{lem-0}).} 
More generally, if $E$ is a twist over $G^Y$, then
\[
\spull E^\phi:=\tphi^*(E)=\set{\(e,\eq[x,z]\)\in E\times
G^X:j(e)=\eq[\phi(x),\phi(z)]}
\]
is a $\T$-groupoid over $G^X$ and we claim that the bundle map
$\ttphi:\spull E^\phi\to E$ given by $\ttphi\(e,\eq[x,z]\)=e$ also induces an
equivalence between $\spull E^\phi$ and $E$.  To verify this, we check
conditions \partref1--\partref5 of \lemref{lem-tech}.  Conditions
\partref1~and \partref2 are straightforward.  For \partref3, suppose
that $\set{\(e_i,\eq[x_i,w_i]\)}$ satisfies
$r\(e_i,\eq[x_i,w_i]\)=\eq[x_i,x_i] \to \eq[x,x]$ and
$s\(e_i,\eq[x_i,w_i]\) = \eq[w_i,w_i] \to \eq[w,w]$ in $G^X$, while
$e_i\to e$ in $E$.  It follows that
$\eq[\phi(x_i), \phi(w_i)]\to \eq[\phi(x),\phi(w)]$ in $G^Y$.  By
passing to a subnet and relabelling, we may assume that there are
$\g_i\in G$ such that
$\(
\phi(x_i)\cdot \g_i,\phi(w_i)\cdot\g_i\)\to \(\phi(x),\phi(w)\)$ in
$Y*Y$.  Similarly, we may assume that there are $\eta_i\in G$ such
that $x_i\cdot \eta_i\to x$ in $X$.  Then $\phi(x_i)\cdot\eta_i\to
\phi(x)$.  Since $Y$ is a principal $G$-space, we can assume that
$\eta_i^{-1}\g_i\to s(x)$.  Thus $x_i\cdot \g_i\to x$; similarly,
$w_i\g_i\to w$ and $\eq[x_i,w_i]\to\eq[x,w]$ in $G^X$; this
establishes~\partref3.

To verify \partref4, suppose that $\eq[x,x]=\eq[w,w]$.  Then
$w=x\cdot\g$ for some $\g\in G$.  But if
$\ttphi\(e,\eq[x,x\cdot\g]\)\in\eo$, then $e\in\eo$ and
$\eq[\phi(x),\phi(x\cdot \g)]=\eq[\phi(x),\phi(x)]$.  Thus $\g=s(x)$
and \partref4 holds.  

Finally if $u=\phi(x)=\phi(w)$, then
$\(i\(\eq[u,u],1\),\eq[x,w]\)\in\ker( \phi)$ and maps
$\(i\(\eq[u,u],1\),\eq[w,w] \)$ to $\(i\([u,u],1\),\eq[x,x]\)$.  Thus
\partref5 is valid and we have proved the claim.

Now suppose that $Y\preceq_\phi X$ and that
$E=\pi^*_{G^Y}(\tb\stt\btb)$ for some principal $\T$-bundle over
$Y/G$.  Then routine computations reveal that $\spull E^\phi\cong
\pi^*_{G^X}(\tb'\stt \btb')$ where $\tb'$ is the pull-back of $\tb$ with
respect to the induced map of $X/G$ onto $Y/G$.  Since $\ttphi$ induces
a homomorphism of $\Tw(G^Y)$ to $\Tw(G^X)$, we obtain a homomorphism
$\Phi$ from $\E(G^Y)$ to $\E(G^X)$.  We want to show that if
$Y\preceq_\psi X$ and if $\Psi$ is the induced homomorphism from
$\E(G^Y)$ to $\E(G^X)$, then $\Psi=\Phi$.  Then we will have a
canonical \hm{} $\Phi_{Y,X}:\E(G^Y)\to\E(G^X)$ whenever $Y\preceq X$.
To show that $\Psi=\Phi$, it will suffice to see that $\spull E^\phi\sim_{G^X}
\spull E^\psi$.  However, we already know, for example, that $\Gr(\ttphi)$ is a
$(\spull E^\phi,E)$-\teq.  Thus $Z:=\Gr(\ttphi)\ste\opposite{\Gr(\ttpsi)}$ is
a $(\spull E^\phi,\spull E^\psi)$-\teq, and we need to see that $Z/\T\cong H$.
However it is straightforward to check that 
\begin{multline*}
Z\cong \set{\(\eq[x,x],
e, \eq[w,w]\)\in X/G\times E\times
X/G: \\
\text{$r(e)=\eq[\phi(x),\phi(w)]$ and
$s(e)=\eq[\phi(w),\phi(w)]$}},
\end{multline*}
 where
\[
(f,\eq[z,x])\cdot\(\eq[x,x],e,\eq[w,w]\) = \(\eq[z,z],fe,\eq[w,w]\)
\]
and 
\[
\(\eq[x,x],e,\eq[w,w]\)\cdot (f,\eq[w,z]\)=
\(\eq[x,x],ef,\eq[z,z]\).
\]
  In particular, 
\begin{multline*}
Z/\T\cong
\set{\(\eq[x,x],h,\eq[w,w]\) \in X/G\times G^Y\times X/G: \\
\text{$r(h)=\eq[\phi(x),\phi(w)]$ and $s(h)=\eq[\phi(w),\phi(w)]$}}.
\end{multline*}
Then $\eq[x,w]\mapsto \(\eq[x,x],\eq[\phi(x),\phi(w)],\eq[w,w]\)$ is
the required isomorphism of $G^X$ and $Z/\T$.

To summarize, we have proved most of the following.

\begin{thm}
\label{thm-extgt}
Suppose that $G$ is a second countable locally compact groupoid. 
If $Y\preceq_\phi X$ in $\pg$, then $\tphi\(\eq[x,w]\):=
\eq[\phi(x),\phi(w)]$
defines a groupoid \hm{} of $G^X$ onto $G^Y$.
If $E$ is a twist over $G^Y$, then let $\spull E^\phi$ be the pull-back by
$\tphi$.
Then $\lbb E \rbb\mapsto \lbb \spull E^\phi
\rbb$ is a well defined \hm{} $\Phi_{Y,X}$ of $\E(G^Y)$ to $\E(G^X)$
which depends only on $X$ and $Y$ and not on the choice of $\phi$.
In particular, $\set{G^X,\Phi_{Y,X}}$ is a directed system of groups.
\end{thm}
\begin{proof}
It only remains to show that if $X\preceq Y \preceq Z$, then
$\Phi_{X,Z}= \Phi_{Y,Z}\circ \Phi_{X,Y}$.  However this follows
immediately from the fact that $\spull{(\spull E^\phi)}^\psi\cong 
\spull E^{\psi\circ \phi}$.
\end{proof}

As described in the introduction, the
two dimensional cohomology of $G$ should be related to twists over $G$
--- or at least to generalized twists over imprimitivity groupoids
$G^X$ for $X\in
\pg$.  Naturally, twists relative to different refinements should be
compared in a common refinement.  Just as in ordinary cohomology
theories, this process is formalized by taking an inductive limit.
This is possible in view of the previous Theorem.

\begin{definition}
If $G$ is a second countable locally compact groupoid, then
\[
\extgt:=\lim_{\substack{\longrightarrow\\X\in\pg}} \E(G^X).
\]
\end{definition}

\else

\fi
\ifdraft
%
%

\section{Third Isomorphism Result}\label{sec-last}

Using our first two isomorphism results (Theorems \ref{thm-main1}~and
\ref{thm-main2}), we see that for each $X\in\pg$, there is an
injective \hm{} $\tau_X:\E(G^X)\to\brg$ given by
\[
\tau_X:= \phi^{\xop}\circ\Th{G^X}^{-1}.
\]
We will write $h_X$ for the canonical \hm s of $\E(G^X)$ into
$\extgt$.

\begin{thm}
\label{thm-main3}
The monomorphisms $\tau_X:\E(G^X)\to\brg$ are compatible with the
directing \hm s $\Phi_{X,Y}:\E(G^X)\to \E(G^Y)$.  In particular, there
is an isomorphism $\tau$ of $\extgt$ onto $\brg$ defined by
\[
\tau\(h_X\(\eclass{E}\)\):=\tau_X\(\eclass E\).
\]
Furthermore, the canonical \hm s $h_X:\E(G^X)\to\extgt$ are injective
for all $X\in\pg$.
\end{thm}

\begin{lem}
\label{lem-0}
If $X\preceq_\phi Y$ and if $Z:=\Gr(\tilde\phi)$, then $Z$ is
isomorphic to $Y\understar G\xop$ as \equi(G^Y,G^X)s.
\end{lem}
\begin{proof}
Each element in $\Gr(\tilde\phi)$ has a representative of the form
$\(\eq[y,y],\eq[\phi(y),x]\)$.  Furthermore the map
$\(\eq[y,y],\eq[\phi(y),x]\) \mapsto \eq[y,x]$ is a well-defined
homeomorphism of $\Gr(\tilde\phi)$ onto $Y\understar G\xop$.  It is
not 
hard to see that this map intertwines the $G^Y$- and $G^X$-actions and
is therefore an isomorphism of equivalences.
\end{proof}

\begin{lem}
\label{lem-2}
If $X\preceq_\phi Y$ and if $Z:=\Gr(\tilde\phi)$, then for all
$[\A,\a]\in\brzg$,
\[
\Phi_{X,Y}\(\Th{G^X}\([\A,\a]\)\)= \Th{G^Y}\(\phi^Z\([\A,\a]\)\).
\]
\end{lem}

\begin{proof}
By \lemref{lem-pi} there is a $\pi\in\rg$ such that
$\aa\sim(X/G\times\K,\bar\pi )$; therefore
$\Th{G^X}\([\A,\a]\)=\eclass{E(\pi)}$ where
$E(\pi)=\set{\(\eq[x,w],U\)\in G^X\times U(\eltwo):
\Ad(U)=\pi_{\smeq[x,w]}}$.  Therefore since $\phi^Z$ is well-defined,
\begin{align*}
\phi^Z\([\A,\a]\) &= [\A^Z,\a^Z] =
\bigl[(X/G\times\K)^Z,{\bar\pi}^Z\bigr] \\
\intertext{which by \lemref{lem-pullb}}
&=\bigl[(\tilde\phi\restr{Y/G})^*(X/G\times\K),\tilde\phi^*\a\bigr]
= [Y/G\times\K,\bar\rho],
\end{align*}
where $\rho_{\smeq[y,z]}=\pi_{\smeq[\phi(y),\phi(z)]}$.
Thus,
\begin{align*}
\Th{G^Y}\circ \phi^Z\([\A,\a]\) &= \eclass{E(\rho)} =
\eclass{\spull{E(\pi)}^\phi} \\
&:= \Phi_{X,Y}\(\eclass{E(\pi)}\) = \Phi_{X,Y}\(\Th{G^X}\([\A,\a]\)\).
\qed
\end{align*}
\renewcommand{\qed}{}
\end{proof}

\begin{proof}
[Proof of \thmref{thm-main3}]
Since the $\tau_X$ are monomorphisms, in order to show that $\tau$ is
a well-defined monomorphism and that the $h_X$ are injective, it will
suffice to see that 
\begin{equation}
\label{eq-***a}
\begin{diagram}
\node{\E(G^X)} \arrow[2]{e,t}{\Phi_{X,Y}}
\arrow{se,b}{\tau_X} \node[2]{\E(G^Y)} \arrow{sw,b}{\tau_Y} \\
\node[2]{\brg}
\end{diagram}
\end{equation}
commutes for all $X\preceq Y$.
But
\begin{align*}
\tau_Y\circ \Phi_{X,Y} 
&= \phi^{\opposite Y}\circ \Th{G^Y}^{-1} \circ \Phi_{X,Y};
\\
\intertext{which, by Lemmas \ref{lem-0}, \ref{lem-2}, and
\ref{lem-inv}\partref2, is equal to}
& \phi^{\opposite Y}\circ \phi^{Y\understar G\xop}\circ \Th{G^X}^{-1}; \\
\intertext{which, by \lemref{lem-inv}\partref5, is equal to}
& \phi^{\xop}\circ\Th{G^X}^{-1} = \tau_X.
\end{align*}
In other words, \eqref{eq-***a} commutes.

It only remains to verify that $\tau$ is surjective.  But if
$[\A,\a]\in\brg$, then there is an $X\in\pg$ such that
$\phi^X\([\A,\a]\) = [\A^X,\a^X]\in \brzg$ (\corref{cor-trivialdd}).
But then \lemref{lem-inv}\partref1\&\partref4 imply that
\[
\tau_X\(\Th{G^X}\([\A^X,\a^X]\)\) = [\A,\a].
\]
This suffices as the image of $\tau$ is the union of the images of the
$\tau_X$ ($X\in\pg$).
\end{proof}
\else

\fi

\ifdraft
%
%

\section{Examples and Additional Remarks}\label{sec-11}

In this section, we collect together several examples that are of interest
and that help to illustrate the general theory we have developed. We also
outline an alternate perspective on our analysis that is revealed by the
third author's analysis of the dual groupoid of a $\cs$-algebra.

We begin by reviewing the examples discussed in the introduction. So,
let $G$ be a locally compact {\em group}.  Since $G$ has but one unit,
the collection of $\cs$-$G$-bundles (Definition~\ref{def-2.14})
coincides with the collection of automorphic representations of $G$
described at the very beginning. Further, since the one point space
consisting of the unit of $G$ carries no cohomology, $\br(G)=\brz (G)$,
and \propref{prop-alex14} is the assertion that $\br(G) $ is really the
collection of exterior equivalence classes of automorphic actions of
$G$. (One does not have to worry about cohomology equivalence of
actions, since for groups, this amounts to unitary equivalence.) Thus
there is complete consistency between the definition of $\br(G)$ that
we gave in the introduction and the definition we gave of $\br(G)$ for
general locally compact groupoids $G$ in the body of this
paper. \thmref{thm-main2}, coupled with the fact that $\br(G)=\brz (G)$
guarantees that $\br(G)\cong \E (G)$.  On the other hand, because
$\go $ reduces to a point, $\E (G)$ reduces to $\Tw(G)$ which,
in turn, is isomorphic to the group of all {\em group\/} extensions of
$\T $ by $G$. Granted that this is isomorphic to $H^2(G,\T)$,
our analysis really does recapture the essence of Mackey {\em et.{}
al}.'s analysis of automorphic representations of locally compact
groups.

However, it does more. The collection of groupoids that are equivalent to a
locally compact group is precisely the collection of all {\em transitive\/}
groupoids. This yields

\begin{prop}
\label{transitive}
If $G$ is a transitive groupoid, then 
\[
\br(G)\cong
\br(G\restr u)\cong H^2(G\restr u,\T ).
\]
 where $G\restr u$ is the isotropy group of
any (and hence every) unit.
\end{prop}

\begin{proof}
This, of course, is a consequence of our Theorem 4.1 and Theorem 2.2B in
\cite{mrw}. (See also Theorem 3.5 in \cite{dakin-seda}.) However, it
should be emphasized that our standing 
assumption about groupoids, that they are locally compact, Hausdorff,
{\em  second countable}, and have {\em Haar systems\/} plays an
important r\^ole here.  Theorem 2.2B of \cite{mrw} says that a transitive
groupoid $G$ satisfying these assumptions is equivalent to the
isotropy group $G\restr u$ of any point $u\in
\go $.  Of course $G\restr u$ is a locally compact (second countable) group.
So, \thmref{thm-main1} and the fact that $\br(G\restr u)\cong H^2(G\restr u,\T )$
complete the proof.
\end{proof}

Turning now to the other ``extreme'', the case of spaces, observe that
a groupoid $G$ is equivalent, in the sense of \cite{mrw}, to a space
precisely when it is principal and acts freely and properly on its
unit space. (Such a groupoid is called a {\em proper principal\/}
groupoid.) Indeed, if $G$ is equivalent to a space $T$, say, then the
fact that $G$ has these properties is straightforward to verify from
the definition of the imprimitivity groupoid of an equivalence
connecting $G$ and $T$.  The converse assertion is Proposition 2.2 of
\cite{mw2}.
 In fact, in this case, $G$ is equivalent to
$\go /G$. Now in the introduction, we already have noted that the
Brauer group of a space $T$, $\br(T)$, is isomorphic to 
$H^2(T,\sT )$ via the Dixmier-Douady map $\delta $. (Here, and throughout this
section we use $\sT $, as before, to denote the sheaf of germs of
continuous $\T $-valued functions on whatever space is appropriate for
the discussion.) This is Green's formulation of Dixmier and Douady's
theory in \cite{green4}. This fact, Proposition 2.2 of \cite{mw2}, and
\thmref{thm-main1} combine to prove

\begin{prop}
\label{principalproper}If $G$ is a proper principal groupoid, then
$\br(G)\cong H^2(\go /G,\sT )$. 
\end{prop}

Several other applications of \thmref{thm-main1} are worth pointing
out. Recall that a groupoid $\Gamma $ is called $r$-discrete if its
range map $r$ is a local homeomorphism. This definition is at slight
odds with \cite{renault}, but now has been universally adopted. The
reason is that $r$ is a local homeomorphism if and only if $\Gamma $
is $r$-discrete in the sense of
\cite{renault} and counting measure is a Haar system for $\Gamma $. In
\cite{alex3},  $r$-discrete groupoids in the sense we are using the
term are called sheaf 
groupoids. In that paper, the first author developed a cohomology theory for 
$r$-discrete groupoids acting on sheaves of abelian groups for the purpose
of identifying Brauer groups. Fix an $r$-discrete groupoid $\Gamma $, let $%
\sT $ denote the sheaf of germs of continuous $\T $-valued
functions on $\Gamma ^{(0)}$, and let $\Gamma $ act on $\sT $ in the
usual way. Then the first author shows how to build a cohomology
theory $H^{*}(\Gamma ,\sT )$, which reduces to Grothendieck's
equivariant cohomology theory \cite{groth}, $H_G^{*}(X,\sT )$, when
$\Gamma $ is the transformation group groupoid $X\times G$ determined
by the action of a discrete group $G$ on a locally compact space
$X$. It may also be worthwhile to point out that certain equivariant
sheaf cohomology groups have been ``computed'' in \cite{groth} (using
spectral sequences) and in \cite{alex3} (using an analog of Milnor's
$\lim^1$ sequence).  Theorem~4.19 of \cite{alex3} may be phrased as

\begin{prop}
If $\Gamma $ is an $r$-discrete groupoid, 
then $\br(\Gamma )\cong H^2(\Gamma,\sT )$.
\end{prop}

Evidently this result contains the case when $\Gamma $ is a space, and so
recovers Green's perspective on the Dixmier-Douady theory. Also, this result
together with our \thmref{thm-main3} suggests that 
for a general locally compact
groupoid $G$ of the kind we have been discussing, there ought to be a
cohomology theory $H^{*}(G,\sT )$ that generalizes $H^{*}(\Gamma
,{\sT})$. In any event, what we can assert is 

\begin{prop}
If $G$ is a locally compact groupoid of the kind we have been discussing,
and if $G$ is equivalent, in the sense of \cite{mrw} to an
$r$-discrete groupoid $
\Gamma $, then $\br(G)\cong H^2(\Gamma ,\sT )$. 
\end{prop}

Of course, the proof is simply to cite \thmref{thm-main1}. However,
the range of 
applicability of this theorem is quite large. Many groupoids of
interest are 
known to be 
equivalent to $r$-discrete groupoids. 
In particular, the holonomy groupoid
of a foliation has this property. See \cite{mrw} for a discussion of
this. Also, 
the groupoids arising in the theory of hyperbolic dynamical systems and
Smale spaces have this property (see \cite{ip96} and \cite{ruelle}).
The exact range of applicability of
this proposition is not known. There are groupoids that are not equivalent
to $r$-discrete groupoids. In fact, a locally compact {\em group\/} is
equivalent to an $r$-discrete groupoid if and only if the group is itself
discrete. (To see this, simply note, as we have noted earlier, that a
groupoid $\Gamma $ is equivalent to a group $G$ iff $\Gamma $ is transitive.
Further, $\Gamma $ is equivalent to any isotropy group $\Gamma \restr
u$. So, if $\Gamma $ is  $r$-discrete, then $G$ is equivalent to the
discrete group $%
\Gamma \restr u$. However, it is easy to see that two groups are
equivalent as 
groupoids if and only if they are isomorphic as topological groups.) 
As far
as we know, the situation for {\em principal\/} groupoids is unknown. There
are, however, shreds of evidence that make the problem quite piquant. For
instance, Ramsay \cite{ramsay2} shows that if $G$ is a locally compact
principal groupoid, then there is a {\em Borel\/} transversal. This means that
there is a Borel set $T\subseteq G^{(0)}$ such that, among other things, the
Borel groupoid $G\restr T$ has countable orbits and is (Borel)
equivalent to $G$ 
(via the Borel equivalence $s^{-1}(T)$). On the other hand,
Schwartzman \cite{schwartzman} gives a necessary and sufficient
condition for a flow to be isomorphic 
to a flow built under a function. (See the discussion on page~282 and note
that the horocyclic flow on the cosphere bundle of a compact Riemann surface
satisfies the condition described there \cite{first}.) Thus there are
flows whose 
associated transformation group groupoids are not equivalent to
${\Z}$-transformation group groupoids via topological transversals;
i.e., a 
continuous flow need not be isomomorphic to a flow built under a function
(cf., \corref{cor-11.6}). 
Whether a transformation group groupoid determined by
a flow is equivalent to an $r$-discrete groupoid in some more complicated
fashion remains to be determined.

\thmref{thm-main1} also yields the main result of
\cite[Theorem~1]{krw}.  Recall, as we
remarked in Remark~\ref{rem-3.4}, if the groupoid under consideration is a
transformation group groupoid, $T\times {\mathfrak G}$, say, then $\br(T\times 
{\mathfrak G})$ is the equivariant Brauer group $\br_{{\mathfrak
G}}(T)$ studied in 
\cite{ckrw}, \cite{krw}, and \cite{prw}.

\begin{prop}[{\cite[Theorem~1]{krw}}]
\label{equivalent}
Suppose that $X$ is a second countable locally compact
Hausdorff space, and that $G$ and $H$ are second countable, locally compact
groups. Suppose that $G$ acts on the left of $X$, while $H$ acts on the
right, so the actions commute, and suppose that both actions are free and
proper. Then the equivariant Brauer groups, $\br_G(X/H)$ and $%
\br_H(G\backslash X)$ are isomorphic.
\end{prop}

\begin{proof}
The point is that the transformation group groupoids, $G\times X/H$ and $%
G\backslash X\times H$, are equivalent by Example 2.4 of \cite{mrw}. Hence,
\thmref{thm-main1} yields the result.
\end{proof}

This result and Corollary 6.1 of \cite{ckrw} yield

\begin{cor}\label{cor-11.6}
Suppose that $X$ is a second countable locally compact space with
countable (integral) cohomology groups, $H^0(X;\Z)$, $H^1(X;\Z)$, and
$H^2(X,\Z)$. If $G$ is the transformation group groupoid associated
with an action of the integers $
\Z $ on $X$, then $\br(G)$ is isomorphic to $H^2(Y,\sT )$ where $Y$
is the space of the flow obtained by suspending the homeomorphism giving the
action of $\Z $ on $X$.
\end{cor}
\begin{remark}
The hypotheses on $X$ are not very restrictive.  For example,
$H^n(X;\Z)$ is countable for all $n$ if $X$ has the homotopy type of a
compact metric space \cite[Lemma~0.3]{rr}. 
\end{remark}

\begin{proof}
Recall how the suspension is obtained. Let $\tau $ be the homeomorphism that
defines the action of $\Z $ on $X$, let $\tilde{Y}=X\times \R $,
and define $\tilde{\tau}$ on $\tilde{Y}$ via the formula $\tilde{\tau}%
(x,t)=(\tau x,t+1)$. Then the $\Z $ action determined by $\tilde{\tau}$
is free and proper and it commutes with the $\R $ action on $\tilde{Y}$
given by translation in the second variable. The space $Y$ is the quotient $%
\tilde{Y}/\Z $. Evidently, the action of $\R $ on $\tilde{Y}$ is
also free and proper. It is clear that $\R \backslash \tilde{Y}{
\times \Z}$ is $X\times \Z $. The suspension of $\tau $ is the
transformation group groupoid $\R \times Y$, that is $\R \times 
\tilde{Y}/\Z $. By Proposition \ref{equivalent}, $\br(G)\cong
\br_{{\R}}(Y)$.  We claim that
our hypotheses on the cohomology of $X$ guarantee
that $H^1(Y;\Z )$ and $H^2(Y;\Z)$
are countable. This will suffice as we can then apply
Corollary 6.1 of \cite{ckrw} to conclude that $\br_{\R }(Y)$ is
isomorphic to $H^2(Y,\sT )$.

To prove the claim, let $q:X\times \R\to Y$ be the quotient map.
Define $Y_1:=q\(X\times(0,\frac34)\)$ and $Y_2:=
q\(X\times(\frac12,\frac54)\)$.  Then $Y_1\cup Y_2=Y$.  Furthermore,
$Y_1$ and $Y_2$ are each homeomorphic to $X\times(0,1)$.
Consequently,  
$H^n(Y_i;\Z)\cong H^n(X;\Z)$ for $i=1,2$; thus the former is
countable provided
$n=1,2$.  Similarly, $Y_1\cap Y_2$ is homeomorphic to a disjoint union
of two copies of $X\times(0,1)$, and $H^n(Y_1\cap Y_2;\Z)$ is
countable for $n=1,2$.  The usual Mayer-Vietoris sequence
\cite[II.5.10]{iversen} gives an exact sequence 
\[
\cdots\to H^{n-1}(Y_1\cap Y_2;\Z) \to H^1(Y_1\cup Y_2;\Z)
\to H^1(Y_1;\Z) \oplus H^1(Y_2;\Z) \to \cdots.
\]
It follows that $H^n(Y_1\cup Y_2;\Z)=H^n(Y;\Z)$ is countable for
$n=1,2$ as claimed.
\end{proof}

Recently there has been considerable interest in the groupoid $G_E$ and
associated groupoid \cs-algebra $\cs(G_E)$ built from the path space
corresponding to a directed graph $E$ \cite{kpr,kprr}.  The
\cs-algebras $\cs(G_E)$ can be quite interesting.  For example, if $A$
is the connectivity matrix of a finite directed graph $E$, then one
obtains the Cuntz-Krieger algebras $\mathcal{O}_A$.

We show that the Brauer group of a groupoid associated to a directed
graph is trivial.  Let $E$ be a (nontrivial) countable directed graph.
We use the notation of \cite{kpr}: the set of vertices is denoted
$E^0$, the set of edges is denoted $E^1$ and the structure maps are
written $r,s : E^1 \to E^0$. We will assume that $E$ is locally finite
(i.e., $r$ and $s$ are finite to one) and that it has no sinks (i.e.,
$s : E^1 \to E^0$ is surjective).  One forms the associated groupoid
$G_E$ as follows: the unit space $\unitspace {G_E}$ is identified with the
infinite path space $E^\infty = \set{ (e_i)_{i=1}^\infty: 
\text{$e_i \in E^1$ and $
r(e_i)=s(e_{i+1})$} }$; the topology of $E^\infty$ is generated by a
countable basis of compact open sets, the cylinder sets associated to
finite paths.  The graph groupoid consists of the collection of
triples
\[
G_E = \set{ (e, n, f) : \text{$e, f \in E^\infty$, $n \in \Z$, and
$e_i = f_{i-n}$ for $i$ sufficiently large} };
\]
the range and source maps are given by $r(e, n, f) = e$ and $s(e, n,
f) = f$, respectively.  Then 
$(e, n, f)(f,m,g) = (e,n+m,g)$ and $(e,n,f)^{-1} = (f,-n,e)$.
The topology of $G_E$ is 
also generated by a countable basis of compact open sets, 
the cylinder sets associated to pairs of 
finite paths ending in a common vertex.  

It will be of some use to have available the following factorization
 of elements of the groupoid which are not units.  For each edge $e_0
 \in E^1$ set $V(e_0) = \set{ (e, 1, f) : \text{$e, f \in E^\infty$, 
$e_1 = e_0$, and $e_i =
 f_{i-1}$ for $i > 1$} } \subset G_E$.  
Given $(e, n, f) \notin
\unitspace {G_E}$, there are $k,l \ge 0$ (at least one of which is nonzero),
 with $n = k - l$ so that $e_{k+i} = f_{l+i}$ for $i \geq 1$ and
 elements $\xi_i \in V(e_i)$ and $\eta_j \in V(f_j)$ for $1 \leq i
 \leq k$ and $
 1 \leq j \leq l$ such that $(e,n,f ) = \xi_1 \dots \xi_k
 {\eta_l}^{-1}\dots {\eta_1}^{-1}$.  There is a unique factorization
 of this form of minimal length.

\begin{prop}
Suppose that $E$ is a locally finite directed graph with no sinks,
and that $G_E$ is the associated locally compact groupoid constructed
above.  Then 
\[
\br(G_E)=\set0.
\]
\end{prop}
\begin{proof}
Since the unit space $E^\infty$ of $G_E$ is zero-dimensional, 
$H^2(E^\infty,{\sT}) = \set0$ 
and so $\br(G_E) = \brz(G_E)$.  
By the second isomorphism result (\thmref{thm-main2}), 
$\brz(G_E)$ is a quotient of $\Tw(G_E)$; thus, it suffices to show
that 
$\Tw(
G_E) = \set0$.  
Suppose that $F$ is a twist over $G_E$ with twist map $j: F \to G_E$;
we will show that there is a trivializing section.
That is, a continuous
section of the twist map which is also a groupoid homomorphism.  Let
$\iota : E^\infty \to F$ be the embedding which identifies the unit
spaces of $G_E$ and $F$.  Since the the twist viewed as a principal
$\T$-bundle is trivial ($G_E$ has a basis of compact open sets),
there is continuous section $\sigma : \bigcup_{x \in E^1} V(x) \to F$ of
the bundle restricted to $\bigcup_{x \in E^1} V(x)$.  We extend $\sigma$
as follows: for $\g \in E^\infty =\unitspace{G_E}$,
 set $\sigma(\g) = \iota(\g)$,
for $\g = (e, n, f) \notin \unitspace{G_E}$, using the above factorization
$(e,n,f) = \xi_1 \dots \xi_k {\eta_l}^{-1}\dots {\eta_1}^{-1}$, set
\[
\sigma(e,n,f) = 
\sigma(\xi_1)\dots\sigma(\xi_k){\sigma(\eta_l)}^{-1}\dots
{\sigma(\eta_1)}^{-1}.
\]
 One checks that $\sigma$ is well-defined, 
and that it is a morphism of groupoids as required.
\end{proof}

We turn now to an alternate viewpoint on the material we have presented that
is based on the concept of the dual groupoid of a $\cs$-algebra. So fix a $%
\cs$-algebra $A$. As a set, the dual groupoid of $A$, denoted $G(A)$,
consists of the extreme points $\omega $ of the unit ball of its dual
$A^{*}$. By taking the left and right polar decompositions of the
functional $%
\omega $ we may write $\omega =xu=uy$, where $x$ and $y$ are pure states and 
$u$ is a partial isometry that is uniquely determined. 
We may therefore
write the elements $\omega \in G(A)$ as triples $\omega =(x,u,y)$. 
Since the partial isometry $u$ in such a
triple is unique, the circle ${\T}$ acts freely on $G(A)$ via the
formula $t(x,u,y):=(x,tu,y)$. This makes $G(A)$ a principal ${\T}$-space
with quotient $R(A)$, the graph of the unitary equivalence relation on the
pure state space of $A$, $P(A)$.
The groupoid structure on $G(A)$ is given by $%
(x,u,y)(y,v,z)=(x,uv,y)$ and $(x,u,y)^{-1}=(y,u^{*},x)$. The
presentation of  
$G(A)$ as triples in this way gives the canonical extension 
\[
P(A)\times \T \arrow{e,J} G(A)\arrow{e,A} R(A).
\]

In the case when $A$ is the \cs-algebra of compact
operators, $\K (\hs)$, on the Hilbert space $\hs$, its pure state space may
be viewed as the projective space $P(\hs):=S(\hs)/\T $, where $S(\hs)$ is the
unit sphere of $\hs$ endowed with the weak topology. The dual groupoid, $%
G(\hs):=G(A)$, then is the quotient of the trivial groupoid $S(\hs)\times S(\hs)$
by the diagonal action of $\T $. It is evident in this case, that the
natural topology on $G(\hs)$ gives $G(\hs)$ the structure of a Polish
groupoid.
The canonical extension
in this case is 
\[
P(\hs)\times \T \arrow{e,J} G(\hs)\arrow{e,A} P(\hs)\times
P(\hs).
\]
We may view the dual groupoid $G(\hs)$ as the groupoid of minimal partial
isometries $u$ with $r(u)=uu^{*}$ and $s(u)=u^{*}u$. 
That is, the map $(x,u,y)\mapsto u$ is an isomorphism of $G(\hs)$
with the collection of all minimal partial isometries on $\hs$.

Let $a\mapsto \gamma (a)$ be an automorphism of a $\cs$-algebra $A$. Then $%
\gamma $ acts on $A^{*}$ by transposition: For $\omega \in A^{*}$, $\omega
\gamma $ is defined by the formula  
$\lip\relax<\omega \gamma ,a>=\lip\relax<\omega ,\gamma (a)>$.  
With this notation, we have $(x,u,y)\gamma =(x\gamma ,\gamma
^{-1}(u),y\gamma )$, for $(x,u,y)\in G(A)$. If $A=\K (\hs)$ and $\alpha $
is a group homomorphism of the locally compact group $G$ into the projective
unitary group $PU(\hs)\cong \Aut(\K (\hs))$, then as we have seen earlier,
the Mackey obstruction $E(\alpha )$ of this projective representation is the
pull back by $\alpha $ of the universal extension 
\[
1\longrightarrow \T \longrightarrow U(\hs)\longrightarrow
PU(\hs)\longrightarrow 1
\]
given by $\Ad$. Explicitly,  
\[
E(\alpha )=\set{(U,\gamma )\in U(H)\times G:\Ad(U)=\alpha (\gamma )}.
\]
The multiplication law is just $(U,\gamma )(V,\gamma ^{\prime })=(UV,\gamma
\gamma ^{\prime })$. As we mentioned in the introduction, this extension
splits if and only if the action is implemented. Indeed, just define
the isomorphism $%
\T \times G\rightarrow E$ by $(z,\gamma )\mapsto (zU(\gamma ),\gamma )$
and conversely. 

The dual groupoid provides an alternative but equivalent definition of the
Mackey obstruction. We view $\alpha $ as an action of $G$ on the \cs%
-algebra $\K (\hs)$. Therefore $G$ acts continuously on $G(\hs)$ on the
right as above. The semi-direct product $G(\hs)*G$ has its multiplication law
given by $(u,\gamma )(\gamma ^{-1}v,\gamma ^{\prime })=(uv,\gamma \gamma
^{\prime })$. One can check that $Z=S(\hs)\times G$ is a $\T $%
-equivalence between $E$ and $G(H)*G$ in the sense of \secref{sec-eh}. Here $E$
acts on the left of $Z$ according to the formula $(U,\gamma )(\xi ,\gamma
^{\prime })=(U\xi ,\gamma \gamma ^{\prime })$ and $G(\hs)*G$ acts on the right
by $(\xi ,\gamma )(\gamma ^{-1}(|\xi \>\<\xi ^{\prime }|),\gamma ^{\prime
})=(\xi ^{\prime },\gamma \gamma ^{\prime })$. (Of course, $|\xi \>\<\xi
^{\prime }|$ denotes the rank one operator that maps $\xi $ to $\xi ^{\prime
}$.) Note that this alternative definition leads us to consider groupoids
and principal $G$-spaces that are no longer locally compact. 

More generally, let $G$ be a locally compact groupoid. Let $({\A},\alpha
)$ be a representative element of $\br(G)$ and write $A=\gotoco(\go ,{\A})$%
. Then the dual groupoid $G(A)$ is a groupoid bundle over $\go $ whose
fiber above $x$ is the dual groupoid $G(A(x))$. The topology on $G(A)$ is
defined naturally in terms of local fields of minimal partial isometries
(which exist by virtue of Fell's condition). For every $\gamma \in G$, the
isomorphism $\alpha _\gamma :A((s(\gamma ))\rightarrow A((r(\gamma ))$
induces the isomorphism $\omega \mapsto \omega \gamma $ from $G(A(r(\gamma
)))$ onto $G(A(s(\gamma )))$ such that $\lip\relax<\omega \gamma ,a>=
\lip\relax<\omega ,\gamma 
(a)>$. If we write as before an element $\omega $ of $G(A(x))$ as a minimal
partial isometry $u\in A(x)$, then $\omega \gamma $ is given by $\gamma
^{-1}(u)$. We define the semi-direct product groupoid $G(A)*G$ as the set of
pairs $(u,\gamma )$ such that $p(u)=r(\gamma )$. The topology is the
relative topology. The multiplication law is given as above: 
\[
(u,\gamma )(\gamma ^{-1}(v),\gamma ^{\prime })=(uv,\gamma \gamma ^{\prime })
\]
and the inverse is given by 
\[
(u,\gamma )^{*}=(\gamma ^{-1}(u),\gamma ^{-1}).
\]
It is a Polish groupoid extension: 
\[
P(A)\times \T \arrow{e,J} G(A)*G\arrow{e,A}
R(A)*G.
\]
Note, in this case, that $R(A)$ is a Polish equivalence relation and that $%
R(A)*G$ is a Polish groupoid equivalent to $G$ via the equivalence $P(A)*G$. 

In our definition of $\extgt $, we consider only principal $G$%
-spaces that are {\em locally compact} and that possess an equivariant $s$%
-system. Therefore, this extension $G(A)*G$ does not quite define an element
of $\extgt $. However, we have the following result:

\begin{prop}
Let $G$ be a locally compact groupoid with Haar system and let $({\A}%
,\alpha )$ be a $\cs$-$G$-system defining an element of $\br(G)$. Then, the
extension $G(A)*G$ defined by the dual groupoid as above is $\T $%
-equivalent, in the sense of \secref{sec-eh}, to any of the extensions provided by 
$({\A},\alpha )$ via the isomorphism of \thmref{thm-main3}.
\end{prop}

\begin{proof}
Let us first assume that $({\A},\alpha )$ determines an element in $%
\brz (G)$. Thus, there exists a continuous field $\HH $ of Hilbert spaces 
$x\mapsto H(x)$ such that $A(x)=\K (H(x))$. The dual groupoid
$G({\HH}):=G(A)$ is then a bundle over $\go $ with fiber
$G(H(x))$. Just as in 
the case of a group discussed above, one can check that $Z:=S(\HH )*G$,
where $S(\HH )$ is the sphere bundle of $\HH $, is a $\T $%
-equivalence between $G(\HH )*G$ and the extension 
\[
E=\{(U,\gamma )\in U(\HH )*G:Ad(U)=\alpha (\gamma )\}
\]
in the sense defined in \secref{sec-eh}. In the general case, we may apply
\corref{cor-trivialdd} to find a groupoid $\underline{G}$ and a
$(\underline{G},G)$%
-equivalence $Z$ such that $({\A}^Z,\alpha ^Z)\in \brz (\underline{G})$.
One can check that $G(A)*Z$ is the desired $\T $-equivalence between $%
G(A)*\underline{G}$ and $G(A^Z)*G$.
\end{proof}

Thus, if one is willing to widen the class of ``equivalent
extensions'' to include groupoids that are not necessarily locally
compact, one finds that the association $({\A},\alpha )\mapsto G(A)*G$
gives a \emph{canonical} representative for an element in $\extgt$
directly in terms of the elementary $\cs$-$G$-bundle
$({\A},\alpha)$. This perspective was emphasized by the third author
in \cite{renault3}, where, in the special case that the groupoid is a 
{\em space}, $T$, he showed how to view the equivalence, $G(A)*Z$, in
the preceding proof as the Dixmier-Douady invariant for
$A=\gotoco(T,{\A})$.

\else

\fi

\includerefs

\end{document}